\documentclass[10pt]{article}

\usepackage[dvips]{epsfig} % this package is for manipulating postscript graphics
\usepackage{amssymb,amsmath,amsthm,mathrsfs} %this package contains all the special ams fonts etc.
\usepackage{graphicx}

\usepackage[authoryear,sort&compress]{natbib}
\usepackage{url}
\usepackage{enumerate}
\usepackage{colordvi,color,pspicture}
\usepackage{psfrag}

\usepackage{setspace}
\graphicspath{{Figures/}}

\newcommand{\Y}{\mathcal{Y}}
\newcommand{\X}{\mathcal{X}}
\newcommand{\U}{\mathcal{U}}
\newcommand{\I}{\mathbf{I}}
\newcommand{\mI}{\mathcal{I}}
\newcommand{\y}{\mathsf{y}}
\newcommand{\NY}{N_{\mathcal{Y}}}
\newcommand{\F}{\mathcal F}
\newcommand{\g}{\gamma}
\newcommand{\G}{\Gamma}
\newcommand{\N}{\mathcal{N}}
\newcommand{\V}{\mathcal{V}}
\newcommand{\A}{\mathcal{A}}
\newcommand{\R}{\mathbb{R}}
\newcommand{\mS}{\mathcal{S}}
\newcommand{\E}{\mathbf{E}}

\newcommand{\D}{\mathscr D}
\newcommand{\ve}{\varepsilon}

\newcommand{\al}{\alpha}
\newcommand{\be}{\beta}

\DeclareMathOperator{\Bernoulli}{Bernoulli}

\theoremstyle{plain}
\newtheorem{theorem}{Theorem}[section]
\newtheorem{lemma}[theorem]{Lemma}
\newtheorem{proposition}[theorem]{Proposition}
\newtheorem{corollary}[theorem]{Corollary}

\theoremstyle{definition}

\newtheorem{example}[theorem]{Example}

\theoremstyle{remark}

\newtheorem{remark}[theorem]{Remark}

\begin{document}
\title{The Distribution of the Domination Number of Class Cover Catch Digraphs for
Non-uniform One-dimensional Data}
\author{Elvan Ceyhan\thanks{Department of Applied Mathematics and Statistics, The Johns Hopkins University,
Baltimore, MD, 21218.}
%\and
%Carey E. Priebe\thanks{Department of Applied Mathematics and Statistics, The Johns Hopkins University,
%Baltimore, Md. 21218}
}
\date{\today}
\maketitle

\maketitle

\begin{abstract}
For two or more classes of points in $\R^d$ with $d \ge 1$,
the class cover catch digraphs (CCCDs)
can be constructed using the relative positions of
the points from one class with respect to the points from the other class.
The CCCDs were introduced by \cite{priebe:2001}
who investigated the case of two classes, $\X$ and $\Y$.
They calculated the exact (finite sample) distribution of the domination number
of the CCCDs based on $\X$ points relative to $\Y$ points both
of which were uniformly distributed on a bounded interval.
We investigate the distribution of the domination number of the CCCDs based
on data from non-uniform $\X$ points on an interval with end points from $\Y$.
Then we extend these calculations for multiple $\Y$ points on bounded intervals.
\end{abstract}

\vspace{0.1 in}

%--------------------------------------------------------------------------------------------------------------------------------------------%\\
\newpage
\section{Introduction}
\label{sec:intro}
In 2001, a new classification method was developed which was based on the relative
positions of the data points from various classes;
\cite{priebe:2001} introduced the class cover catch digraphs (CCCDs)
in $\R$ and gave the exact distribution of the domination number of the CCCDs
for two classes, $\X$ and $\Y$, with uniform distribution on a bounded interval in $\R$.
\cite{devinney:2002b} proved a SLLN result for the one-dimensional
class cover problem.
\cite{devinney:2002a}, \cite{marchette:2003}, \cite{priebe:2003b}, and \cite{priebe:2003a}
extended the CCCDs to higher dimensions and demonstrated that CCCDs
are a competitive alternative to the existing methods in classification.
The classification method based on CCCDs involves data reduction (condensing)
by using approximate --- rather than exact --- minimum dominating sets as \emph{prototype sets},
since finding the exact minimum dominating set for CCCDs is an NP-hard problem in general.
However for finding a dominating set of CCCDs on the real line,
a simple linear time algorithm is available (\cite{priebe:2001}).
But unfortunately, the exact and the asymptotic distributions of the domination number
of the CCCDs are not analytically tractable in multiple dimensions.

To address the latter issue of intractability of the distribution
of the domination number in multiple dimensions,
\cite{ceyhan:2003a,ceyhan:2005e} introduced the central similarity proximity maps
and $r$-factor proportional-edge proximity maps and the associated random
proximity catch digraphs.
Proximity catch digraphs are a generalization of the CCCDs.
The asymptotic distribution of the domination number of the latter is calculated and then used in
testing spatial patterns between two or more classes.
See \cite{ceyhan:2005e} for more detail.

In this article, we generalize the original result of \cite{priebe:2001}
to the case of non-uniform $\X$ points with support being the interval with end points from $\Y$,
and then to multiple $\Y$ points in a bounded interval $(c,d) \subset \R$ with $c<d$.
These generalizations will also serve as the bases for extension of the
results for the uniform and non-uniform data in higher dimensions.

\section{Data-random Class Cover Catch Digraphs}
Let $(\Omega,\mathcal M)$ be a measurable space and $\X_n=\{X_1,\ldots,X_n\}$ and
$\Y_m=\{Y_1,\ldots,Y_m\}$ be two sets of $\Omega$-valued random variables
from classes $\X$ and $\Y$, respectively, with joint probability distribution $F_{X,Y}$.
Let $d(\cdot,\cdot):\Omega\times \Omega \rightarrow [0,\infty)$ be any distance function.
The \emph{class cover problem} for a target class, say $\X$, refers to finding a collection of neighborhoods,
$N_i$ around $X_i$ such that
(i) $\X_n \subseteq  \bigl(\cup_i N_i \bigr)$ and (ii) $\Y_m \cap \bigl(\cup_i N_i \bigr)=\emptyset$.
A collection of neighborhoods satisfying both conditions is called a {\em class cover}.
A cover satisfying condition (i) is a {\em proper cover} of class $\X$
while a cover satisfying condition (ii) is a {\em pure cover} relative to class $\Y$.
This article is on the {\em minimum cardinality class covers}; that is,
class covers satisfying both (i) and (ii) with the smallest number of neighborhoods.
See \cite{priebe:2001}.

Consider the map $N:\Omega \rightarrow 2^\Omega$
where $2^\Omega$ represents the power set of $\Omega$.
Then given $\Y_m \subseteq \Omega$,
the {\em proximity map}
$\NY(\cdot): \Omega \rightarrow 2^\Omega$
associates with each point $x \in \Omega$
a {\em proximity region} $\NY(x) \subseteq \Omega$.
For $B \subseteq \Omega$, the $\G_1$-region is the image of the map
$\G_1(\cdot,\NY):2^\Omega \rightarrow 2^\Omega$
that associates the region $\G_1(B,\NY):=\{z \in \Omega: B \subseteq  \NY(z)\}$
with the set $B$.
For a point $x \in \Omega$, we denote $\G_1(\{x\},\NY)$ as $\G_1(x,\NY)$.
Notice that while the proximity regions are defined for one point,
$\G_1$-regions are defined for sets of points.

The {\em data-random CCCD} has the vertex set $\V=\X_n$
and arc set $\A$ defined by $(X_i,X_j) \in \A \iff X_j \in \NY(X_i)$.
In particular, we use $\NY(X_i)=B(X_i,r_i)$, the open ball around $X_i$ with radius
$r_i:=\min_{Y \in \Y_m}d(X_i,Y)$, as the proximity map as in \cite{priebe:2001}.
We call such a digraph a $\D_{n,m}$-digraph.
A $\D_{n,m}$-digraph is a {\em pseudo digraph} according some authors
if loops are allowed (see, e.g., \cite{chartrand:1996}).

A data-random CCCD for $\Omega=\R^d$ and $N_i=B(X_i,r_i)$
is referred to as $\mathscr C_{n,m}$-graph in \cite{priebe:2001}.
We change the notation to emphasize the fact that $\D_{n,m}$ is a digraph.
Furthermore, \cite{ceyhan:2003a} call the proximity map $N_i=B(X_i,r_i)$
a {\em spherical proximity map}.

The $\D_{n,m}$-digraphs are closely related to the {\em proximity graphs} of
\cite{jaromczyk:1992} and might be considered as a special case of
{\em covering sets} of \cite{tuza:1994} and {\em intersection digraphs} of \cite{sen:1989}.
Our data-random proximity digraph is a {\em vertex-random proximity digraph} and
not a standard one (see e.g., \cite{janson:2000}).
The randomness of a $\D_{n,m}$-digraph lies in the fact that
the vertices are random with the joint distribution $F_{X,Y}$,
but arcs $(X_i,X_j)$ are
deterministic functions of the random variable $X_j$ and the random set $N_i$.

\section{Domination Number of Random $\D_{n,m}$-digraphs}
In a digraph $D=(\V,\A)$ of order $|\V|=n$, a vertex $v$ {\em dominates}
itself and all vertices of the form $\{u:\,(v,u) \in \A\}$.
A {\em dominating set}, $S_D$, for the digraph $D$ is a subset of
$\V$ such that each vertex $v \in \V$ is dominated by a vertex in $S_D$.
A {\em minimum dominating set},  $S^*_D$, is a dominating set of minimum cardinality;
and the {\em domination number}, denoted $\g(D)$, is defined as $\g(D):=|S^*_D|$,
where $|\cdot|$ is the set cardinality functional (\cite{west:2001}).
If a minimum dominating set consists of only one vertex,
we call that vertex a {\em dominating vertex}.
The vertex set $\V$ itself is always a dominating set,
so $\g(D) \le n$.

Let $\F\left( \R^d \right):=\{F_{X,Y} \text{ on } \R^d \text { with } P(X=Y)=0\}$.
As in \cite{priebe:2001}, in this article,
we consider $\D_{n,m}$-digraphs for which
$\X_n$ and $\Y_m$ are random samples from $F_X$ and $F_Y$, respectively,
and the joint distribution of $X,Y$ is $F_{X,Y} \in \F\left( \R^d \right)$.
We call such digraphs as \emph{$\F\left( \R^d \right)$-random $\D_{n,m}$-digraphs}
and focus on the random variable $\g(D)$.
To make the dependence on sample sizes explicit,
we use $\g(D_{n,m})$ instead of $\g(D)$.
It is trivial to see that $1 \le \g(D_{n,m}) \le n$,
and $\g(D_{n,m}) < n$ for nontrivial digraphs.

\section{The Distribution of the Domination Number of $\F(\R)$-random $\D_{n,m}$-digraphs}
\label{sec:domination-number-Dnm}
In $\R$, the data-random CCCD is a special case of
{\em interval catch digraphs} (see, e.g., \cite{sen:1989} and \cite{prisner:1994}).
Let $\X_n$ and $\Y_m$ be two samples from $\F(\R)$ and
$Y_{(j)}$ be the $j^{th}$ order statistic of $\Y_m$ for $j=1,2,\ldots,m$.
Then $Y_{(j)}$ partition $\R$ into $(m+1)$ intervals.
Let
$$-\infty =: Y_{(0)}<Y_{(1)}< \ldots <Y_{(m)}< Y_{(m+1)}:=\infty,$$
and $\mI_j:=\left(Y_{(j-1)},Y_{(j)}\right)$, $\X^j:=\X_n \cap \mI_j$,
and $\Y^j:=\{Y_{(j-1)},Y_{(j)}\}$ for $j=1,2,\ldots,(m+1)$.
This yields a disconnected digraph with subdigraphs $D^j$ for $j=1,2,\ldots,(m+1)$,
each of which might be null or itself disconnected.
Let $\g(D^j)$ denote the the cardinality of the minimum dominating set
for the component of the random $\D_{n,m}$-digraph
induced by the pair $\X^j$ and $\Y^j$, $n_j:=|\X^j|$,
and $F_j$ be the density $F_X$ restricted to $\mI_j$.
Then $\g(D_{n,m})=\sum_{j=1}^{m+1}\g(D^j)$.
We study the simpler random variable $\g(D^j)$ first.
The following lemma follows trivially (see \cite{priebe:2001}).

\begin{lemma}
\label{lem:end-intervals}
For $j\in\{ 1,(m+1)\}$, $\g(D^j)=\I(n_j >0)$ where $\I(\cdot)$ is the indicator function.
\end{lemma}

For $j=2,\ldots,m$ and $n_j > 0$, we prove that $\g(D^j) \in \{1,2\}$
with the distribution dependent probabilities $1-p_{n_j}(F_j),p_{n_j}(F_j)$, respectively,
where $p_{n_j}(F_j)=P(\g(D^j)=2)$.
A quick investigation shows that $\g(D^j)=2$ iff
$\X^j \cap \left( \frac{\max\,(\X^j)+Y_{(j-1)}}{2},\frac{\min(\X^j)+Y_{(j)}}{2} \right) =\emptyset$;
that is, $\X^j \subset B(x,r(x))$ iff
$x \in \left( \frac{\max\,(\X^j)+Y_{(j-1)}}{2},\frac{\min(\X^j)+Y_{(j)}}{2} \right)$
where $r(x)=\min(x-Y_{(j-1)},Y_{(j)}-x)$.
Hence $\G_{1}(\X^j,\NY)=\left( \frac{\max\,(\X^j)+Y_{(j-1)}}{2},\frac{\min(\X^j)+Y_{(j)}}{2} \right)\subseteq \mI_j$.
By definition, if $\X^j \cap \G_1(\X^j,\NY) \not=\emptyset$,
then $\g(D^j)=1$; hence the name \emph{$\G_1$-region} and the notation $\G_{1}(\cdot,\NY)$.

\begin{theorem}
\label{thm:gamma 1 or 2}
For $j=2,\ldots,m$, $\g(D^j) \sim 1+\Bernoulli(p_{n_j}(F_j))$ for $n_j > 0$.
\end{theorem}

{\bf Proof:}
See \cite{priebe:2001} for the proof.
$\blacksquare$

The probability $P(\g(D^j)=2)=P(\X^j \cap \G_{1}(\X^j,\NY) =\emptyset)$
depends on the conditional distribution $F_{X|Y}$ and the interval $\G_1(\X^j,\NY)$,
which, if known, will make possible the calculation of $p_{n_j}(F_j)$.
As an immediate result of Lemma \ref{lem:end-intervals} and Theorem \ref{thm:gamma 1 or 2},
we have the following upper bound for $\g(D_{n,m})$.

\begin{theorem}
Let $D_{n,m}$ be an $\F(\R)$-random $\D_{n,m}$-digraph with $n>0,\,m>0$
and $k_1$ and $k_2$ be two natural numbers defined as
$k_1:=\sum_{j=2}^m \I(|\X_n \cap \mI_j|>1)$
and
$k_2:=\sum_{j=2}^m \I(|\X_n \cap \mI_j|=1)+\sum_{j \in \{1,(m+1)\}}
\I(\X_n \cap \mI_j\not=\emptyset)$.
Then $1\le \g(D_{n,m})\le 2\,k_1+k_2 \le \min(n,2\,m)$.
\end{theorem}

In the special case of fixed $\Y_2=\{\y_1,\y_2\}$ and $\X_n$
a random sample from $\U(\y_1,\y_2)$, the uniform distribution on $(\y_1,\y_2)$,
we have a $\D_{n,2}$-digraph
for which $F_X=\U(\y_1,\y_2)$ and $F_Y$ is a degenerate distribution.
%with an atom at singleton $(\y_1,\y_2)$.
We call such digraphs as \emph{$\U(\y_1,\y_2)$-random $\D_{n,2}$-digraphs}
and provide an exact result on the distribution of their domination
number in the next section.

\subsection{The Exact Distribution of the Domination Number of $\U(\y_1,\y_2)$-random
$\D_{n,2}$-digraphs}
\label{sec:gamma-dist-uniform}
Suppose $\Y_2=\{\y_1,\y_2\} \subset \R$ with $-\infty<\y_1<\y_2<\infty$ and
$\X_n =\{X_1,\ldots,X_n\}$ a set of iid random variables from $\U(\y_1,\y_2)$.
%(e.g., $X_i$ are from a Poisson point process).
Any $\U(\y_1,\y_2)$ random variable can be transformed into a $\U(0,1)$
random variable by $\phi(x)=(x-\y_1)/(\y_2-\y_1)$,
which maps intervals $(t_1,t_2) \subseteq (\y_1,\y_2)$ to
intervals $\bigl( \phi(t_1),\phi(t_2) \bigr) \subseteq (0,1)$.
So, without loss of generality, we can assume $\X_n=\{X_1,\ldots,X_n\}$
is a set of iid random variables from the $\U(0,1)$ distribution.
That is, the distribution of $\g(D_{n,2})$ does not depend on the support interval $(\y_1,\y_2)$.
%Let $\g_n:=\g(D_{n,2})$ for brevity of notation.
Recall that $\g(D_{n,2})=2$ iff $\X_n \cap \G_1(\X_n,\NY)=\emptyset$,
then $P(\g(D_{n,2})=2)=4/9-(16/9) \, 4^{-n}$.
For more detail, see (\cite{priebe:2001}).
Hence, for $\U(\y_1,\y_2)$ data, we have
\begin{equation}
\label{eqn:finite-sample-unif}
\g(D_{n,2})=  \left\lbrace \begin{array}{ll}
       1           & \text{w.p. $5/9+(16/9) \, 4^{-n}, $}\\
       2           & \text{w.p. $4/9-(16/9) \, 4^{-n}, $}
\end{array} \text{~~for all $n \ge 1$}, \right.
\end{equation}
where w.p. stands for ``with probability".
Then the asymptotic distribution of $\g(D_{n,2})$ for $\U(\y_1,\y_2)$ data is given by
\begin{equation}
\label{eqn:asymptotic-uniform}
\lim_{n \rightarrow \infty}\g(D_{n,2}) =
\left\lbrace \begin{array}{ll}
       1           & \text{w.p. $5/9$,}\\
       2           & \text{w.p. $4/9$.}
\end{array} \right.
\end{equation}
For $m>2$, \cite{priebe:2001} computed the exact distribution of $\g(D_{n,m})$.
However, independence of the distribution of the domination number from the support interval
does not hold in general;
that is, for $X_i \stackrel{iid}{\sim}F$ with support $\mS(F) \subseteq(\y_1,\y_2)$,
the exact and asymptotic distribution of $\g(D_{n,2})$ will
depend on $F$ and $\Y_2$.

\subsection{The Distribution of the Domination Number for $\F(\R)$-random $\D_{n,2}$-digraphs}
\label{sec:non-uniform}

For $\Y_2=\{\y_1,\y_2\} \subset \R$ with $-\infty<\y_1<\y_2<\infty$,
a quick investigation shows that the $\G_1$-region is
$\G_1(\X_n,\NY)=\left( \frac{\y_1+X_{(n)}}{2},\frac{\y_2+X_{(1)}}{2} \right)$.
Note that $\X_n \cap \G_1(\X_n,\NY)$ is the set of all dominating vertices,
which is empty when $\g(D_{n,2})>1$.
To make the dependence on $F$ explicit and for brevity of notation,
we will denote the domination number of the
$F\bigl((\y_1,\y_2)\bigr)$-random $\D_{n,2}$-digraphs as $\g_n(F)$.

Let $p_n(F):=P(\g_n(F)=2)$ and $p(F):=\lim_{n \rightarrow \infty}P(\g_n(F)=2)$.
Then the exact (finite sample) and asymptotic distributions of $\g_n(F)$
are $1+\Bernoulli\left(p_n(F)\right)$ and $1+\Bernoulli\left(p(F)\right)$, respectively.
That is, for finite $n$, we have
\begin{equation}
%\label{eqn:finite-sample-unif}
\g_n(F)=  \left\lbrace \begin{array}{ll}
       1           & \text{w.p. $1-p_n(F)$}\\
       2           & \text{w.p. $p_n(F)$}
\end{array} \text{for all $n \ge 1$}. \right.
\end{equation}
The asymptotic distribution is similar.

With $\Y_2=\{0,1\}$, let $F$ be a distribution with support $\mS(F) \subseteq (0,1)$ and density $f$
and let $\X_n$ be a set of $n$ iid random variables from $F$.
Since $\g_n(F) \in \{1,\,2\}$, to find the distribution of $\g_n(F)$,
it suffices to find $P(\g_n(F)=1)$ or $P\bigl( \g_n(F)=2 \bigr)$.
For computational convenience, we employ the latter in our calculations.

Then
{\small
\begin{equation}
\label{eqn:prob-g=2}
p_n(F)=\int_{\mS(F) \setminus \G_1(\X_n,\NY)}
\left[ 1-\frac{ F((1+x_1)/2)-F(x_n/2)}{F(x_n)-F(x_1)}\right]^{n-2}f_{1n}(x_1,x_n)dx_ndx_1,
\end{equation}
}
where
$f_{1n}(x_1,x_n)=n\,(n-1)\,\bigl[ F(x_n)-F(x_1) \bigr]^{n-2}\,f(x_1)\,f(x_n)\;\I(0 < x_1<x_n<1)$
which is the joint probability density function of $X_{(1)},X_{(n)}$.

If the support $\mS(F)=(0,1)$, then the region of integration becomes
$$\Bigl\{ (x_1,x_n) \in (0,1)^2:\; (1+x_1)/2 \le x_n \le 1;\;\;
0\le x_1 \le 1/3 \text{ or }  2\,x_1 \le x_n \le 1;\;\;1/3 \le x_1 \le 1/2  \Bigr\}.$$
The integrand in Equation \eqref{eqn:prob-g=2} simplifies to
\begin{equation}
\label{eqn:integrand}
H(x_1,x_n):=n\,(n-1)f(x_1)f(x_n)\bigl[F(x_n)+F\left(x_n/2 \right)-
\left( F\left( (1+x_1)/2 \right)+F(x_1) \right)\bigr]^{n-2}.
\end{equation}

Let $\X_n$ be a set of iid random variables from a continuous distribution
$F$ with $\mS(F)\subseteq (0,1)$.
The simplest of such distributions is $\U(0,1)$, the uniform distribution
on $(0,1)$, which yields the simplest exact distribution for $\g_n(F)$.
If $X \sim F$, then by probability integral transform, $F(X) \sim \U(0,1)$.
So for any continuous $F$,
we can construct a proximity map depending on $F$ for which
the distribution of the domination number for the
associated digraph will have the same distribution
as that of $\g_n(\U(0,1))$.
%In this case, we design proximity maps that yield the same distribution as $\g_n(\U(0,1))$.

\begin{proposition}
Let $X_i \stackrel{iid}{\sim} F$ which is an (absolutely) continuous
distribution with support $\mS(F)=(0,1)$ and $\X_n=\{X_1,\ldots,X_n\}$.
Define the proximity map $N_F(x):=F^{-1}(\NY(F(x)))=F^{-1}(B(F(x),r(F(x))))$
where $r(F(x))=\min(F(x),1-F(x))$.
Then the domination number of the digraph based on $N_F$, $\X_n$, and $\Y_2=\{0,1\}$,
is equal in distribution to $\g_n(\U(0,1))$.
\end{proposition}
{\bfseries Proof:}
Let $U_i:=F(X_i)$ for $i=1,\ldots,n$ and $\U_n=\{U_1,\ldots,U_n\}$.
Hence, by probability integral transform, $U_i \stackrel{iid}{\sim} \U(0,1)$.
Let $U_{(k)}$ be the $k^{th}$ order statistic of $\U_n$ for $k=1,\ldots,n$.
Furthermore, such an $F$ preserves order; that is, for $x \le y$, $F(x) \le F(y)$.
So the image of $N_F(x)$ under $F$ is
$F(N_F(x))=\NY(F(x))=B(F(x),r(F(x)))$ for (almost) all $x \in (0,1)$.
Then $F(N_F(X_i))=\NY(F(X_i))=\NY(U_i)$ for $i =1,\ldots,n$.
Since $U_i \stackrel{iid}{\sim} \U(0,1)$, the distribution of
the domination number of the digraph based on $\NY$, $\U_n$ and $\{0,1\}$
is given in Equation \eqref{eqn:finite-sample-unif}.
Observe that
$X_j \in N_F(X_i)$ iff
$X_j \in F^{-1}(B(F(X_i),r(F(X_i))))$ iff
$F(X_j) \in B(F(X_i),r(F(X_i)))$ iff
$U_j \in B(U_i,r(U_i))$ for $i,j=1,\ldots,n$.
Hence $P(\X_n \subset N_F(X_i))=P(\U_n \subset \NY(U_i))$
for all $i=1,\ldots,n$.
Therefore, $\X_n \cap \G_1(\X_n,N_F)=\emptyset$ iff $\U_n \cap \G_1(\U_n,\NY)=\emptyset$,
which implies that the domination number of the digraph based on
$N_F$, $\X_n$, and $\Y_2=\{0,1\}$ is 2 with probability $4/9-(16/9)\,4^{-n}$.
Hence the desired result follows.
$\blacksquare$

For example for $F(x)=x^2\,\I(0\le x \le 1)+ \I(x>1)$,
\begin{equation*}
N_F(x)=  \left\lbrace \begin{array}{ll}
       \left( 0,\sqrt{2}\,x \right)      & \text{for $x \in \left[0,1/\sqrt{2}\,\right],$}\vspace{.025in}\\
       \left( \sqrt{2\,x^2-1},1 \right)  & \text{for $x \in \left(1/\sqrt{2},1\,\right].$}
\end{array}  \right.
\end{equation*}

There is also a stochastic ordering between $\g_n(F)$ and $\g_n(\U(0,1))$
provided that $F$ satisfies some conditions
which are given in the following proposition.

\begin{proposition}
Suppose $\X_n=\{X_1,\ldots,X_n\}$ is a random sample from
a continuous distribution $F$ with $\mS(F)\subseteq(0,1)$ and let $X_{(j)}$
be the $j^{th}$ order statistic of $\X_n$ for $j=1,\ldots,n$.
If
\begin{equation}
\label{eqn:stoch-order}
F\bigl( X_{(n)}/2 \bigr)< F\left(X_{(n)}\right)/2 \text{ and }
F\bigl( X_{(1)} \bigr) < 2\,F\left( \left( 1+X_{(1)} \right)/2 \right)-1 \text{ hold a.s., }
\end{equation}
then $\g_n(F)<^{ST}\g_n(\U(0,1))$.
If $<$'s in expression \eqref{eqn:stoch-order} are replaced with $>$'s,
then $\g_n(F)>^{ST}\g_n(\U(0,1))$.
If $<$'s in expression \eqref{eqn:stoch-order} are replaced with $=$'s,
then $\g_n(F)\stackrel{d}{=}\g_n(\U(0,1))$ where $\stackrel{d}{=}$ stands
for equality in distribution.
\end{proposition}
\noindent
{\bfseries Proof:}
Let $U_i:=F(X_i)$ for $i=1,\ldots,n$ and $\U_n=\{U_1,\ldots,U_n\}$.
Then, by probability integral transform, $U_i \stackrel {iid}{\sim} \U(0,1)$.
Let $U_{(j)}$ be the $j^{th}$ order statistic of $\U_n$ for $j=1,\ldots,n$.
The $\G_1$-region for $\U_n$ based on $\NY$ is
$\G_1(\U_n,\NY)=\left(U_{(n)}/2,\,\left(1+U_{(1)}\right)/2 \right)$;
likewise, $\G_1(\X_n,\NY)=\left(X_{(n)}/2,\,\left(1+X_{(1)}\right)/2 \right)$.
But the conditions in expression \eqref{eqn:stoch-order} imply that
$\G_1(\U_n,\NY) \subsetneq F(\G_1(\X_n,\NY))$.
So $\U_n \cap F(\G_1(\X_n,\NY)) = \emptyset$ implies that
$\U_n \cap \G_1(\U_n,\NY) = \emptyset$ and
$\U_n \cap F(\G_1(\X_n,\NY)) = \emptyset$ iff $\X_n \cap \G_1(\X_n,\NY) = \emptyset$.
Hence
$$p_n(F)=P(\X_n \cap \G_1(\X_n,\NY) = \emptyset)<
P(\U_n \cap \G_1(\U_n,\NY) = \emptyset)=p_n(\U(0,1)).$$
Then $\g_n(F)<^{ST}\g_n(\U(0,1))$ follows.
The other cases can be shown similarly.
$\blacksquare$

For more on the comparison of $\g_n(F)$ for general $F$ against $\g_n(\U(0,1))$,
see Section 4.2.2 of the technical report by \cite{ceyhan:2004c}.

\subsubsection{The Exact Distribution of $\g_n(F)$ for $F$ with
Piecewise Constant Density}
\label{sec:piecewise-constant}
Let $\Y_2=\{0,1\}$.
We can find the exact distribution of $\g_n(F)$ for $F$ whose density is piecewise constant.
Note that the simplest of such distributions is the uniform distribution $\U(0,1)$.
Below we give some examples for such densities.
%Denote $p_n(F)$ as $p_n(F)$ for brevity of notation.

\begin{example}
\label{ex:PW1}
Consider the distribution $F$ with density $f(\cdot)$ which is of the form
$f(x)=\frac{1}{1-2\,\delta}\,\I\bigl( \delta < x < 1-\delta \bigr)
\text{ with }\delta \in [0,1/2).$
%\begin{figure} [ht]
%    \centering
%   \scalebox{.29}{\input{Piecewise1.pstex_t}}
%   \scalebox{.29}{\input{Piecewise2.pstex_t}}
%    \caption{The graphs of the density in Example \ref{ex:PW1} (left) and Example \ref{ex:PW2} (right).}
%\label{fig:PW1-2}
%\end{figure}
%See Figure \ref{fig:PW1-2} (left).
Then $F(x)=\frac{x-\delta}{1-2\,\delta}\,\I\bigl( \delta < x < 1-\delta \bigr)+
\I\bigl( x\ge 1-\delta \bigr).$
The integrand in Equation \eqref{eqn:integrand} becomes
$$ H(x_1,x_n)=\frac{n(n-1)}{(1-2\,\delta)^2}
\left(\frac{3\,(x_n-x_1)-1}{2\,(1-2\,\delta)}\right)^{n-2}.$$
Then for $\delta \in [0,1/3]$
\begin{eqnarray}
\label{eqn:pwc1}
p_n(F)&=&\int_{\delta}^{1/3}\int_{(1+x_1)/2}^{1-\delta} H(x_1,x_n)\,dx_ndx_1+
\int_{1/3}^{(1-\delta)/2}\int_{2\,x_1}^{1-\delta}H(x_1,x_n) \,dx_ndx_1 \nonumber\\
 &=& \left(4/9-(16/9)\,4^{-n}\right)\left( \frac{1-3\,\delta}{1-2\,\delta} \right)^n,
\end{eqnarray}
which converges to $0$ as $n \rightarrow \infty$ at (an exponential) rate
$O\bigl( (\frac{1-3\,\delta}{1-2\,\delta})^n \bigr)$.
For $\delta \in [1/3,1/2)$,  it is easy to see that $\g_n(F)=1$ a.s.
In fact, for $\delta \in [1/3,1/2)$ the corresponding digraph is a complete digraph of
order $n$, since $\X_n \subset N(X_i)$ for each $i=1,\ldots,n$.
Furthermore, if $\delta=0$,
then $F=\U(0,1)$ which yields $p_n(F)=4/9-(16/9)\,4^{-n}$.
$\square$
\end{example}

\begin{example}
\label{ex:PW2}
Consider the distribution $F$ with density $f(\cdot)$ which is of the form
$$f(x)=\frac{1}{1-2\,\delta}\,\I\bigl(x \in (0,1) \setminus (1/2-\delta,1/2+\delta)
\bigr) \text{ with }\delta \in [0,1/6].$$
%See Figure \ref{fig:PW1-2} (right).
Then the cumulative distribution function (cdf) is given by
$$F(x)=F_1(x)\,\I\bigl(0<x<1/2-\delta\bigr)+F_2\,(x)\,
\I\bigl(1/2-\delta<x<1/2+\delta\bigr)+F_3\,(x)\,\I\bigl(1/2+\delta<x<1 \bigr)+\I\bigl( x \ge 1 \bigr),$$
where
$$F_1(x)=x/(1-2\,\delta),\;\;F_2\,(x)=1/2,\;\; \text{ and } F_3\,(x)=
(x-2\,\delta)/(1-2\,\delta).$$
There are four cases regarding the relative position of
$X_{(n)}/2, \bigl( 1+X_{(1)} \bigr)/2$ and
$1/2-\delta,1/2+\delta$ that yield $\g_n(F)=2$:
\begin{align*}
\text{{\bf case (1)}}&
\left( X_{(n)}/2, \bigl( 1+X_{(1)} \bigr)/2 \right) \subseteq \bigl(1/2-\delta,1/2+\delta \bigr);&
\text{{\bf case (2)}}&~ X_{(n)}/2<1/2-\delta< \bigl( 1+X_{(1)} \bigr)/2 <1/2+\delta;\\
\text{{\bf case (3)}}&~ 1/2-\delta<X_{(n)}/2<1/2+\delta< \bigl( 1+X_{(1)} \bigr)/2;&
\text{{\bf case (4)}}&~ X_{(n)}/2<1/2-\delta<1/2+\delta< \bigl( 1+X_{(1)} \bigr)/2.
\end{align*}
Let $E_j(n)$ be the event for which \textbf{case (j)} holds for $j=1,2,3,4$,
for example,
$$E_1(n):=\left\{\left( X_{(n)}/2, \bigl( 1+X_{(1)} \bigr)/2 \right) \subseteq
\bigl(1/2-\delta,1/2+\delta \bigr) \right\}.$$
Then $p_n(F)=\sum_{j=1}^4 P\left( \g_n(F)=2,\,E_j(n) \right)$.
Furthermore, \textbf{cases (2)} and \textbf{(3)} are symmetric; i.e.,
$P(\g_n(F)=2,E_2(n))=P(\g_n(F)=2,E_3(n))$.
Then in {\bf case (1)},
we obtain $P(\g_n(F)=2,E_1(n))=1-2\,\left(\frac{1-4\,\delta}{1-2\,\delta}\right)^n+
\left(\frac{1-6\,\delta}{1-2\,\delta}\right)^n$.
%{\bf case (1)}:\\
%If $\left( X_{(n)}/2, \bigl( 1+X_{(1)} \bigr)/2 \right)
%\subseteq \bigl(1/2-\delta,1/2+\delta \bigr)$, i.e.,
%$\G_1(\X_n,\NY) \subseteq \bigl( 1/2-\delta,1/2+\delta \bigr)$,
%then $\X_n \cap \G_1(\X_n,\NY)=\emptyset$, which implies $\g_n(F)=2$.
%Hence $P(\g_n(F)=2,\,E_1(n))=P(E_1(n))$
%and $P(E_1(n)) \le P\bigl( \g_n(F)=2 \bigr)$.
%Then for $\delta \in (0,1/6)$
%$$P\bigl(E_1(n)\bigr)=
%\int_{0}^{2\,\delta}\int_{1-2\,\delta}^1 f_{1n}(x_1,x_n)\,dx_ndx_1 =
%1-2\,\left(\frac{1-4\,\delta}{1-2\,\delta}\right)^n+\left(\frac{1-6\,\delta}{1-2\,\delta}\right)^n$$
%where
%\begin{eqnarray*}
%f_{1n}(x_1,x_n)&=&n\,(n-1)\,\bigl( F_3\,(x_n)-F_1(x_1) \bigr)^{n-2}\,f(x_1)\,f(x_n)\\
%&=&n\,(n-1)\,\bigl( x_n-x_1-2\,\delta \bigr)^{n-2}\,\bigl( 1-2\,\delta \bigr)^{-n}.
%\end{eqnarray*}
Note that $P\bigl(\G_1(\X_n,\NY) \subseteq \bigl( 1/2-\delta,1/2+\delta \bigr)\bigr)\rightarrow 1$
as $n \rightarrow \infty$, hence it suffices to use this case to show that
$p_n(F) \rightarrow 1$ as $n \rightarrow \infty$ at an exponential rate
since $P(E_1(n))\le p_n(F)$.

In \textbf{cases (2)} and \textbf{(3)}, we obtain
$P(\g_n(F)=2,E_2(n))=\frac{2}{3}\, \left(1-\frac{4}{4^n} \right)\,
\left( \left(\frac{1-4\,\delta}{1-2\,\delta}\right)^n-
\left(\frac{1-6\,\delta}{1-2\,\delta}\right)^n \right)$ and in
\textbf{case (4)}, $P(\g_n(F)=2,E_4(n))=\frac{4}{9}\,\bigl(
1-4^{-n+1} \bigr)\,\left(\frac{1-6\,\delta}{1-2\,\delta}\right)^n.$
See \cite{ceyhan:2004c} for the details of the computations.

Combining the results from the cases, for $\delta \in [0,1/6]$ we have
\begin{equation}
\label{eqn:pwc2}
P\bigl( \g_n(F)=2 \bigr)=1+\left(\frac{1-6\,\delta}{1-2\,\delta}\right)^n\,
\bigl(1/9+(32/9)4^{-n}\bigr)-\left(\frac{1-4\,\delta}{1-2\,\delta}\right)^n\,
\bigl(2/3+(16/3)4^{-n} \bigr),
\end{equation}
which converges to 1 as $n \rightarrow \infty$ at rate
$O\left( \left(\frac{1-4\,\delta}{1-2\,\delta}\right)^n \right)$.

Notice that if $\delta=0$, then $F=\U(0,1)$.
The exact distribution for $\delta \in (1/6,1/3)$ can be found in a similar fashion.
Furthermore, if $\delta \in [1/3,1/2]$, then $p_n(F)=1-2\,\delta^n$.
See \cite{ceyhan:2004c} also for the details of the computations.
$\square$
\end{example}

\begin{example}
\label{ex:PW3}
Consider the distribution $F$ with density $f(\cdot)$ which is of the form
$f(x)=(1+\delta)\,\I\bigl(x \in (0,1/2)\bigr) +(1-\delta)\,\I\bigl(x \in [1/2,1)\bigr) \text{ with }\delta \in [-1,1].$

%\begin{figure} [ht]
%   \centering
%   \scalebox{.3}{\input{Piecewise3a.pstex_t}}
%   \scalebox{.3}{\input{Piecewise4b.pstex_t}}
%   \caption{The plot of the density in Example \ref{ex:PW3} with $\delta>0$ (left) and Example \ref{ex:PW4} with $\delta>0$ %(right).}
%\label{fig:PW3-4}
%\end{figure}

%See Figure \ref{fig:PW3-4} (left) for $f(x)$ with a positive $\delta$ value.
Then
\begin{equation}
\label{eqn:pwc3}
p_n(F)=\frac{4(1-\delta^2)}{9-\delta^2}-\frac{8\cdot 4^{-n}(1-\delta^2)}{3}
\left(\frac{(1+\delta)^{n-1}}{3-\delta}+\frac{(1-\delta)^{n-1}}{3+\delta}\right).
\end{equation}
See \cite{ceyhan:2004c} for the derivation.
Hence $ \lim_{n \rightarrow \infty} p_n(F)= \frac{4\,(1-\delta^2)}{9-\delta^2} =:p_F(\delta),$
with the rate of convergence $O\left( \left(\frac{1+\delta}{4}\right)^n  \right)$.
Note that $p_F(\delta) \in [0,4/9]$ is continuous in $\delta$ and decreases as $|\delta|$ increases.
If $\delta=0$, then $F=\U(0,1)$ and $p_F(\delta=0)=4/9$.
Note also that $p_F(\delta=\pm 1)=0$.
$\square$
\end{example}

\begin{example}
\label{ex:PW4}
Consider the distribution $F$ with density $f(\cdot)$ which is of the form
$$f(x)=(1+\delta)\,\I(0<x<1/4)+(1-\delta)\,\I(1/4 \le x< 3/4)+
(1+\delta)\,\I(3/4 \le x< 1) \text{ with $\delta \in [-1,1]$}.$$
%See Figure \ref{fig:PW3-4} (right) for $f(x)$ with a positive $\delta$ value.

The exact value of $p_n(F)$ is available,
but it is rather a lengthy expression (see \cite{ceyhan:2004c} for the expression
and its derivation).
But the limit is as follows:
$p_n(F)\rightarrow \frac {4\,(1+\delta)^2}{(3+\delta)^2}=:p_F(\delta)$
as $n \rightarrow \infty$ with the rate of convergence
$O\left( \left(\frac{5-\delta}{8}\right)^n  \right)$.
So $p_F(\delta)$ is increasing in $\delta$.
Notice here that $p_n(F)$ and $p_F(\delta)$ are continuous in $\delta$ and
$p_F(\delta)>0$ for all $\delta \in (-1,1]$.
Moreover, $p_F(\delta=1)=1$ and $p_F(\delta=-1)=0$.
$\square$
\end{example}

Note that extra care should be taken if the points of discontinuity in the above examples
are different from $\{1/4,\,3/4\}$ or $1/2$, since the symmetry in the
probability calculations no longer exists in such cases.

\subsubsection{The Exact Distribution of $\g_n(F)$ for Polynomial $f$ Using Multinomial Expansions}
\label{sec:polynomial}
The exact distribution of $\g_n(F)$ for (piecewise) polynomial
$f(x)$ with at least one piece is of degree 1 or higher can be
obtained using the multinomial expansion of the term $(\cdot)^{n-2}$
in Equation \eqref{eqn:integrand} with careful bookkeeping.
However, the resulting expression for $p_n(F)$ is extremely lengthy
and not that informative.

The simplest example is with $f(x)=2\,x$ and $F(x)=x^2$.
Then $p_n(F)=P\bigl( \g_n(F)=2 \bigr)=\Lambda_1(n)+\Lambda_2(n),$
where
$\Lambda_1(n):=\int_{0}^{1/3}\int_{(1+x_1)/2}^{1}H(x_1,x_n)dx_ndx_1$,
$\Lambda_2(n):=\int_{1/3}^{1/2}\int_{2\,x_1}^{1}H(x_1,x_n)dx_ndx_1$,
and
$H(x_1,x_n)=n\,(n-1) {x_1}\,{x_n}\, \bigl(5\,x_n^2-1-2\,x_1-5\,x_1^2 \bigr)^{n-2}.$
Then
\begin{equation*}
\Lambda_1(n)=\int_0^{1/3}\bigl( 8\,n\,x_1/5 \bigr)\,\bigl(1-x_1/2-5\,x_1^2/4 \bigr)^{n-1}
-\bigl( 8\,n\,x_1/5 \bigr)\, \bigl( 1/16+x_1/2-15\,x_1^2/16 \bigr)^{n-1}\,dx_1.
\end{equation*}
Using the multinomial expansion of $(\cdot)^{n-1}$ with respect to $x_1$
in the integral above,
we have
$$\Lambda_1(n)=\sum_{Q_2}{n-1 \choose q_1,q_2,q_3}\frac{8\,n\,(-1)^{q_2+q_1}{5}^{-1+q_1}{2}^{-q_2-2\,q_1}{3}^{-2-q_2-2\,q_1}}{2+q_2+2\,q_1}+\frac{n\,(-1)^{1+q_1}{2}^{3-3\,q_2-4\,q_3-4\,q_1}{15}^{q_1}{3}^{-2-q_2-2\,q_1}}{5\,(2+q_2+2\,q_1)}$$
where $Q_2=\bigl\{q_1,q_2,q_3 \in {\mathbb N}: q_1+q_2+q_3=n-1 \bigr\}$.

Similarly, the second piece follows as
\begin{equation*}
\Lambda_2\,(n)=\int_{1/3}^1\bigl( 8\,n\,x_1/5 \bigr)\,\bigl( 1-x_1/2-5\,x_1^2/4 \bigr)^{n-1}
-\bigl( 8\,n\,x_1/5 \bigr)\, \bigl(15/x_1^2/4-1/4-x_1/2 \bigr)^{n-1}\,dx_1.
\end{equation*}

Again, using the multinomial expansion of the $(\cdot)^{n-1}$ term above, we get
\begin{multline*}
\Lambda_2\,(n)=\sum_{Q_3} {n-1 \choose r_1,r_2,r_3}\Bigl[2\,n \Bigl(
9\,(-1)^{{r_2}+r_1}5^{r_1}4^{-2\,r_1-r_2}+
9\,(-1)^{1+{r_3}+r_2}15^{r_1}4^{-2\,r_1-r_3-r_2}+\\
4\,(-1)^{1+r_2+r_1}6^{-r_2-2\,r_1}5^{r_1}+(-1)^{r_3+r_2}{4}^{1-r_3}
6^{-r_2}12^{-r_1}5^{r_1}\Bigr)\Bigr]\Big/\Bigl[90+45\,r_2+90\,r_1 \Bigr]
\end{multline*}
where $Q_3=\bigl\{ r_1,r_2,r_3 \in {\mathbb N}: r_1+r_2+r_3=n-1 \bigr\}$.
See \cite{ceyhan:2004c} for more detail and examples.

%It is also possible to obtain $p_n(F)$ by numerical integration
For fixed numeric $n$, one can obtain $p_n(F)$ for $F$ (omitted for the sake of brevity)
with the above densities by numerical integration of the below expression.
\begin{eqnarray*}
p_n(F)=P\bigl( \g_n(F)=2 \bigr)&=&\int_{0}^{1/3}\int_{(1+x_1)/2}^{1}
H(x_1,x_n)+\int_{1/3}^{1/2}\int_{2\,x_1}^{1}H(x_1,x_n)\,dx_ndx_1,
\end{eqnarray*}
where $H(x_1,x_n)$ is given in Equation $\eqref{eqn:integrand}$.

Recall the $\F(\R^d)$-random $\D_{n,m}$-digraphs.
We call the digraph which obtains
in the special case of $\Y_m=\{\y_1,\y_2\}$ and support of $F_X$ in $(\y_1,\y_2)$,
\emph{$\F((\y_1,\y_2))$-random $\D_{n,2}$-digraph}.
Below, we provide asymptotic results
pertaining to the distribution of such digraphs.

\section{The Asymptotic Distribution of the Domination Number of
$\F((\y_1,\y_2))$-random $\D_{n,2}$-digraphs}
Although the exact distribution of $\g_n(F)$ is not analytically available
in a simple closed form for $F$
whose density is not piecewise constant,
the asymptotic distribution of $\g_n(F)$ is available for larger families of distributions.
%by resorting to asymptotics in probability calculations.
% Equation \eqref{eqn:prob-g=2}.
 First, we present the asymptotic distribution of $\g_n(F)$ for
$\D_{n,2}$-digraphs with $\Y_2=\{\y_1,\y_2\} \subset \R$ with $\y_1<\y_2$
for various $F$ with support $\mS(F) \subseteq (\y_1,\y_2)$.
Then we will extend this to the case with $\Y_m \subset \R$ for $m>2$.
% and $\mS(F) \subseteq (\y_1,\y_2)$ in Section \ref{sec:general-support}.

For $\ve \in (0,(\y_1+\y_2)/2)$, consider the family of distributions given by
\begin{multline*}
\F\bigl((\y_1,\y_2),\ve \bigr)=\Bigl \{\text{$F$ : %$\mS(F)=(\y_1,\y_2)$ and
$(\y_1,\y_1+\ve) \cup (\y_2-\ve,\y_2)\cup \bigl( (\y_1+\y_2)/2-\ve,(\y_1+\y_2)/2+\ve \bigr) \subseteq \mS(F)\subseteq(\y_1,\y_2)$} \Bigr\}.
\end{multline*}

Let the $k^{th}$ order right (directed) derivative at $x$ be defined as
$f^{(k)}(x^+):=\lim_{h \rightarrow 0^+}\frac{f^{(k-1)}(x+h)-f^{(k-1)}(x)}{h}$
for all $k \ge 1$ and the right limit at $c$ be defined as $f(c^+):=\lim_{h \rightarrow 0^+}f(c+h)$.
The left derivatives and limits are defined similarly with $+$'s being replaced by $-$'s.
Furthermore,
let $\vec{h}=(h_1,h_2)$ and $\vec{c}=(c_1,c_2)$ and
the directional limit at $(c_1,c_2) \in \R^2$ for $g(x,y)$ in the first quadrant in $\R^2$ be
$g(c_1^+,c_2^+):=\lim_{\substack{||\vec{h}|| \rightarrow 0\\
h_1,h_2>0}}g(\vec{c}+\vec{h})$
and the directional partial derivatives at $(c_1,c_2)$ along paths in the first quadrant be
$\frac{\partial^{k+1} g(c_1^+,c_2^+)}{\partial x^{k+1}} := \lim_{\substack{||\vec{h}|| \rightarrow 0\\
h_1,h_2>0}}
\frac{1}{||h||}\left( \frac{\partial^k g(\vec{c}+h)}{\partial x^k}-\frac{\partial^k g(\vec{c})}{\partial x^k} \right)$
for $k\ge 1$.

\begin{theorem}
\label{thm:kth-order-gen}
Let $\Y_2=\{\y_1,\y_2\} \subset \R$ with $-\infty < \y_1 < \y_2<\infty$
and $\X_n=\{X_1,\ldots,X_n\}$ with $X_i \stackrel {iid}{\sim} F \in \F((\y_1,\y_2),\ve)$ .
Let $D_{n,2}$ be the random $\D_{n,2}$-digraph based on $\X_n$ and $\Y_2$.
Suppose $k \ge 0$ is the smallest integer for which
$F(\cdot)$ has continuous right derivatives up to order $(k+1)$ at $\y_1,\,(\y_1+\y_2)/2$,
$f^{(k)}(\y_1^+)+2^{-(k+1)}\,f^{(k)}\left( \left( \frac{\y_1+\y_2}{2} \right)^+ \right) \not= 0$
and $f^{(j)}(\y_1^+)=0$ for all $j=0,1,\ldots,k-1$;
and $\ell \ge 0$ is the smallest integer for which
$F(\cdot)$ has continuous left derivatives up to order $(\ell+1)$ at $\y_2,\,(\y_1+\y_2)/2$,
$f^{(\ell)}(\y_2^-)+2^{-(\ell+1)}\,f^{(\ell)}\left( \left( \frac{\y_1+\y_2}{2} \right)^- \right) \not= 0$
and $f^{(j)}(\y_2^-)=0$ for all $j=0,1,\ldots,\ell-1$.
Then $\g_n(F) \sim 1+ \Bernoulli\bigl(p_n(F)\bigr)$ where
$p_n(F):=P\bigl(\g_n(F)=2\bigr)$ and for bounded $f^{(k)}(\cdot)$ and $f^{(\ell)}(\cdot)$,
we have the following limit
$$\lim_{n \rightarrow \infty}p_n(F) =
\frac{f^{(k)}(\y_1^+)\,f^{(\ell)}(\y_2^-)}{\left[f^{(k)}(\y_1^+)+2^{-(k+1)}\,f^{(k)}
\left( \left( \frac{\y_1+\y_2}{2} \right)^+ \right)\right]\,\left[f^{(\ell)}(\y_2^-)+
2^{-(\ell+1)}\,f^{(\ell)}\left( \left( \frac{\y_1+\y_2}{2} \right)^- \right)\right]}.$$
Note also that $p_1(F)=0$.
\end{theorem}

{\bfseries Proof:}
First suppose $(\y_1,\y_2)=(0,1)$.
Recall that $\G_1(\X_n,\NY)=\left(X_{(n)}/2,
\bigl( 1+X_{(1)} \bigr)/2 \right) \subset (0,1)$
and $\g_n(F)=2 \text{ iff } \X_n \cap \G_1(\X_n,\NY)=\emptyset$.
Then for finite $n$,
\begin{equation*}
p_n(F)=P\bigl( \g_n(F)=2 \bigr)=\int_{\mS(F) \setminus \G_1(\X_n,\NY)} H(x_1,x_n)\,dx_ndx_1,
\end{equation*}
where $H(x_1,x_n)$ is as in Equation \eqref{eqn:integrand}.
%Above the integrand simplifies to
%$$n\,(n-1)f(x_1)f(x_n)\left[F(x_n)-F(x_1)+F\left(x_n/2 \right)-F\left( (1+x_1)/2 \right)\right]^{n-2}.$$

Let $\ve \in (0,1/3)$.
Then $P\bigl( X_{(1)} <\ve,\; X_{(n)} > 1-\ve \bigr)\rightarrow 1$
as $n \rightarrow \infty$ with the rate of convergence depending on $F$.
So for sufficiently large $n$,
\begin{equation}
\label{eqn:asy-g=2-x1xn}
p_n(F) \approx \int_0^{\ve}\int_{1-\ve}^1 n\,(n-1)f(x_1)f(x_n)
\Bigl[F(x_n)-F(x_1)+F\left(x_n/2 \right)-F\left( (1+x_1)/2 \right) \Bigr]^{n-2}\,dx_ndx_1.
\end{equation}
 Let
$$G(x_1,x_n)=F(x_n)-F(x_1)+F\left(x_n/2 \right)-F\left( (1+x_1)/2 \right).$$
The integral in Equation \eqref{eqn:asy-g=2-x1xn} is critical at $(x_1,x_n)=(0,1)$,
since $G(0,1)=1$
and for $(x_1,x_n) \in (0,1)^2$ the
integral converges to 0 as $n \rightarrow \infty$.
So we make the change of variables $z_1=x_1$ and $z_n=1-x_n$, then $G(x_1,x_n)$ becomes
$$G(z_1,z_n)=F(1-z_n)-F(z_1)+F\left( (1-z_n)/2 \right)-F\left( (1+z_1)/2 \right),$$
and Equation \eqref{eqn:asy-g=2-x1xn} becomes
\begin{equation}
\label{eqn:asy-g=2-z1zn}
p_n(F) \approx \int_0^{\ve}\int_0^{\ve} n\,(n-1)f(z_1)f(1-z_n)
\left[G(z_1,z_n)\right]^{n-2}\,dz_ndz_1.
\end{equation}
The new integral is critical at $(z_1,z_n)=(0,0)$. Note that
$\frac{\partial^{r+s}G(z_1,z_n)}{\partial z_1^r \, \partial z_n
^s}=0$ for all $r,s \ge 1$. Let $\al_i := \frac{\partial^{i+1}
G(z_1,z_n)}{\partial z_1^{i+1}}|_{(0^+,0^+)}=
f^{(i)}(0^+)+2^{-(i+1)}\,f^{(i)}\left(\frac{1}{2}^+\right)$ and
$\be_j := \frac{\partial^{j+1} G(z_1,z_n)}{\partial
z_n^{j+1}}|_{(0^+,0^+)}=
f^{(j)}(1^-)+2^{-(j+1)}\,f^{(j)}\left(\frac{1}{2}^-\right)$. Then by
the hypothesis of the theorem, we have $\al_i = 0$ and
$f^{(i)}\left(\frac{1}{2}^+\right)=0$ for all $i=0,1,\ldots,(k-1)$;
and $\be_j = 0$ and $f^{(j)}\left(\frac{1}{2}^-\right)=0$ for all
$j=0,1,\ldots,(\ell-1)$.
So the Taylor series expansions of $f(z_1)$
around $z_1=0^+$ up to order $k$ and $f(1-z_n)$ around $z_n=0^+$ up
to order $\ell$, and $G(z_1,z_n)$ around $(0^+,0^+)$ up to order
$(k+1)$ and $(\ell+1)$ in $z_1,z_n$, respectively, so that
$(z_1,z_n) \in (0,\ve)^2$, are as follows. {\small
\begin{multline*}
f(z_1)=\frac{1}{k!}f^{(k)}(0^+)\,z_1^k+O\left( z_1^{k+1} \right);~~~~
f(1-z_n)=\frac{(-1)^\ell}{\ell!}f^{(\ell)}(1^-)\,z_n^\ell+O\left(z_n^{\ell+1}\right);\\
G(z_1,z_n)=G(0^+,0^+)+ \frac{1}{(k+1)!}\left(\frac{\partial^{k+1}
G(0^+,0^+)}{\partial z_1^{k+1}}
\right)\,z_1^{k+1}+\frac{1}{(\ell+1)!}\left(\frac{\partial^{\ell+1} G(0^+,0^+)}
{\partial z_n^{\ell+1}}\right)\,z_n^{\ell+1}
+O\left( z_1^{k+2} \right) + O\left(z_n^{\ell+2}\right)\\
=1-\frac{\al_k}{(k+1)!}\,z_1^{k+1}
+\frac{(-1)^{\ell+1}\be_\ell}{(\ell+1)!}\,z_n^{\ell+1}+O\left(z_1^{k+2}\right) +
O\left(z_n^{\ell+2}\right).
\end{multline*}
}
%where $G(0^+,0^+)=\lim_{h \rightarrow 0^+}G(h,h)$ and
%$\frac{\partial^{k+1} G(0^+,0^+)}{\partial z_1^{k+1}}= \lim_{h \rightarrow 0^+}
%\frac{1}{h}\left( \frac{\partial^{k} G(h,h)}{\partial z_1^{k}}-\frac{\partial^{k} G(0^+,0^+)}{\partial z_1^{k}} \right)$.
Then substituting these expansions in Equation \eqref{eqn:asy-g=2-z1zn}, we obtain
\begin{multline*}
p_n(F) \approx \int_0^{\ve}\int_0^{\ve}
n(n-1)\Biggl[\frac{1}{k!}f^{(k)}(0^+)\,z_1^k+O\left(z_1^{k+1}\right)\Biggr]
\Biggl[\frac{(-1)^\ell}{\ell!}f^{(\ell)}(1^-)\,z_n^\ell +O\left(z_n^{\ell+1}\right)\Biggr]\\
\Biggl[1-\frac{\al_k}{(k+1)!}\,z_1^{k+1}
-\frac{(-1)^{\ell}\be_\ell}{(\ell+1)!}\,z_n^{\ell+1}+O\left(z_1^{k+2}\right) +
O\left(z_n^{\ell+2}\right)\Biggr]^{n-2}\,dz_ndz_1.
\end{multline*}

Now we let $z_1=w\,n^{-1/(k+1)}$, $z_n=v\,n^{-1/(\ell+1)}$,
%$m=(k+\ell+1)\,\min(1/(k+1),1/(\ell+1))$,
and $\nu=\min(k,\ell)$ to obtain
%{\small
\begin{multline}
\label{eqn:Pg2-in-(0,1)} p_n(F) \approx
\int_0^{\ve\,n^{1/(k+1)}}\int_0^{\ve\,n^{1/(\ell+1)}}n\,(n-1)
\Biggl[\frac{1}{n^{k/(k+1)}\,k!}f^{(k)}(0^+)w^k+O\left(n^{-1}\right)\Biggr]
\Biggl[\frac{(-1)^\ell}{n^{\ell/(\ell+1)}\,\ell!}f^{(\ell)}(1^-)v^\ell
+O\left(n^{-1}\right)\Biggr]\\
\Biggl[1-\frac{1}{n}\left(\frac{\al_k}{(k+1)!}\,w^{k+1}+
\frac{(-1)^{\ell}\be_\ell}{(\ell+1)!}\,v^{\ell+1}\right)+ O\left(
n^{-(\nu+2)/(\nu+1)}\right)\Biggr]^{n-2}\,\left(\frac{1}{n^{1/(k+1)}}\right)\,
\left(\frac{1}{n^{1/(\ell+1)}}\right)\,dvdw\\
 =  \int_0^{\ve\,n^{1/(k+1)}}\int_0^{\ve\,n^{1/(\ell+1)}} n\,(n-1)
\Biggl[\frac{(-1)^\ell}{n^2\,k!\,\ell!} f^{(k)}(0^+)
f^{(\ell)}(1^-)w^k v^\ell+O\left(n^{-(2k+3)/(k+1)}\right)+
O\left(n^{-(2\ell+3)/(\ell+1)}\right)\\
+O\left(n^{-2(k+2)(\ell+2)/((k+1)(\ell+1))}\right)\Biggr]
\Biggl[1-\frac{1}{n}\Bigl[\frac{\al_k}{(k+1)!}\,w^{k+1}+
\frac{(-1)^{\ell}\,\be_\ell}{(\ell+1)!}\,v^{\ell+1}\Bigr]+
O\left( n^{-(\nu+2)/(\nu+1)} \right)\Biggr]^{n-2}\,dvdw,\\
\text{letting $n \rightarrow
\infty,$~~~~~~~~~~~~~~~~~~~~~~~~~~~~~~~~~~~~~~~~~~~~~~~~~~~~~~~~
~~~~~~~~~~~~~~~~~~~~~~~~~~~~~~~~~~~~~~~~~~~~~~~~~~~~}\\
\approx  \int_0^{\infty}\int_0^{\infty} \frac{(-1)^\ell}{k!\,\ell!}
f^{(k)}(0^+) f^{(\ell)}(1^-)w^k v^\ell \,
\exp\Biggl[-\frac{\al_k}{(k+1)!}\,w^{k+1}-
\frac{(-1)^{\ell}\,\be_\ell}{(\ell+1)!}\,v^{\ell+1}\Biggr]\,dvdw\\
= \frac{f^{(k)}(0^+)\, f^{(\ell)}(1^-)\,(-1)^\ell\,(k+1)!(\ell+1)!}
{k!\,\ell!\,(-1)^\ell\,(k+1)(\ell+1)\,\al_k\,\be_\ell}
= \frac{f^{(k)}(0^+)\, f^{(\ell)}(1^-)}{\al_k\,\be_\ell}\\
= \frac{f^{(k)}(0^+)\,f^{(\ell)}(1^-)}{\left[f^{(k)}(0^+)+
2^{-(k+1)}\,f^{(k)} \left( \frac{1}{2}^+ \right)\right]\,
\left[f^{(\ell)}(1^-)+2^{-(\ell+1)}\,f^{(\ell)} \left( \frac{1}{2}^-
\right)\right]},
\end{multline}
%}
as $n \rightarrow \infty$ at rate $O(c(f)\cdot n^{-m})$ where $c(f)$
is a constant depending on $f$.

For the general case of $\Y=\{\y_1,\y_2\}$, the transformation
$\phi(x)=\frac{x-\y_1}{\y_2-\y_1}$ maps $(\y_1,\y_2)$ to $(0,1)$ and
the transformed random variables $\phi(X_i)$ are distributed with
density $g(x)=(\y_2-\y_1)\,f\left( \frac{x-\y_1}{\y_2-\y_1} \right)$
on $(0,1)$. Substituting $f(x)$ by $g(x)$ in Equation
\eqref{eqn:Pg2-in-(0,1)}, the desired result follows. $\blacksquare$

Note that
\begin{itemize}
\item
if $\min\bigl( f^{(k)}(\y_1^+),f^{(\ell)}(\y_2^-) \bigr)=0$ and
$\min\left( f^{(k)} \left( \frac{(\y_1+\y_2)}{2}^+ \right),
\,f^{(\ell)} \left( \frac{(\y_1+\y_2)}{2}^- \right) \right)\not=0$
then $p_n(F)\rightarrow 0$ as $n \rightarrow \infty$, at rate $O\bigl( c(f)\cdot n^{-m} \bigr)$ and
\item
if $\min\bigl( f^{(k)}(\y_1^+),f^{(\ell)}(\y_2^-) \bigr)\not=0$ and
$f^{(k)} \left( \frac{(\y_1+\y_2)}{2}^+ \right)=
f^{(\ell)} \left( \frac{(\y_1+\y_2)}{2}^- \right)=0$
then $p_n(F) \rightarrow 1$
as $n \rightarrow \infty$, at rate $O\bigl(c(f)\cdot n^{-m} \bigr)$.
\end{itemize}
For example, with $F=\U(\y_1,\y_2)$, in Theorem \ref{thm:kth-order-gen} we have
$k=\ell=0$, $f(\y_1^+)=f(\y_2^-)=f\left( \frac{(\y_1+\y_2)}{2}^+ \right)=
f\left( \frac{(\y_1+\y_2)}{2}^- \right)=1/(\y_2-\y_1)$.
Then $\lim _{n \rightarrow \infty}p_n(F)=4/9$,
which agrees with the result given in Equation \eqref{eqn:asymptotic-uniform}.

\begin{example}
For $F$ with density $f(x)=\bigl( x+1/2 \bigr)\,\I\bigl( 0 <x<1 \bigr)$,
we have $k=\ell=0$, $f(0^+)=1/2$, $f(1^-)=3/2$ and $f\left( \frac{1}{2}^+ \right)=
f\left( \frac{1}{2}^- \right)=1$.
Thus $\lim_{n\rightarrow \infty}p_n(F)=3/8 =0.375$.
The numerically computed (by numerical integration) value of $p_n(F)$ with
$n=1000$ is $\widehat p_{1000}(F)\approx 0.3753$.
$\square$
\end{example}

\begin{remark}
Let $p_F:=\lim_{n \rightarrow \infty}p_n(F)$.
Then the finite sample mean and variance of $\g_n(F)$ are given by
$1+p_n(F)$ and $p_n(F)\,(1-p_n(F))$, respectively;
and the asymptotic mean and variance of $\g_n(F)$ are given by
$1+p_F$ and $p_F\,(1-p_F)$, respectively. $\square$
\end{remark}

\begin{remark}
\label{rem:unbounded}
In Theorem \ref{thm:kth-order-gen}, we assume that $f^{(k)}(\cdot)$ and $f^{(\ell)}(\cdot)$
are bounded on $(\y_1,\y_2)$.
Suppose either $f^{(k)}(\cdot)$ or $f^{(\ell)}(\cdot)$ or both are
not bounded on $(\y_1,\y_2)$ for $k,l \ge 0$,
in particular at $\y_1,(\y_1+\y_2)/2,\y_2$, for example,
$\lim_{x \rightarrow \y_1^+}f^{(k)}(x)=\infty$.
Then we find $p(F)$ as
$$p(F) =\lim_{\delta \rightarrow 0^+}\frac{f^{(k)}(\y_1+\delta)\,f^{(\ell)}(\y_2-\delta)}{\left[f^{(k)}(\y_1+\delta)+2^{-(k+1)}\,f^{(k)} \left( \frac{(\y_1+\y_2)}{2}+\delta \right)\right]\,\left[f^{(\ell)}(\y_2-\delta)+2^{-(\ell+1)}\,f^{(\ell)} \left( \frac{(\y_1+\y_2)}{2}-\delta \right)\right]}.\;\;\square$$
\end{remark}

\begin{figure}[ht]
\centering
\psfrag{x}{\scriptsize{$x$}}
\psfrag{f(x)}{\scriptsize{$f(x)$}}
\psfrag{density #}{}
\rotatebox{-90}{\epsfig{figure=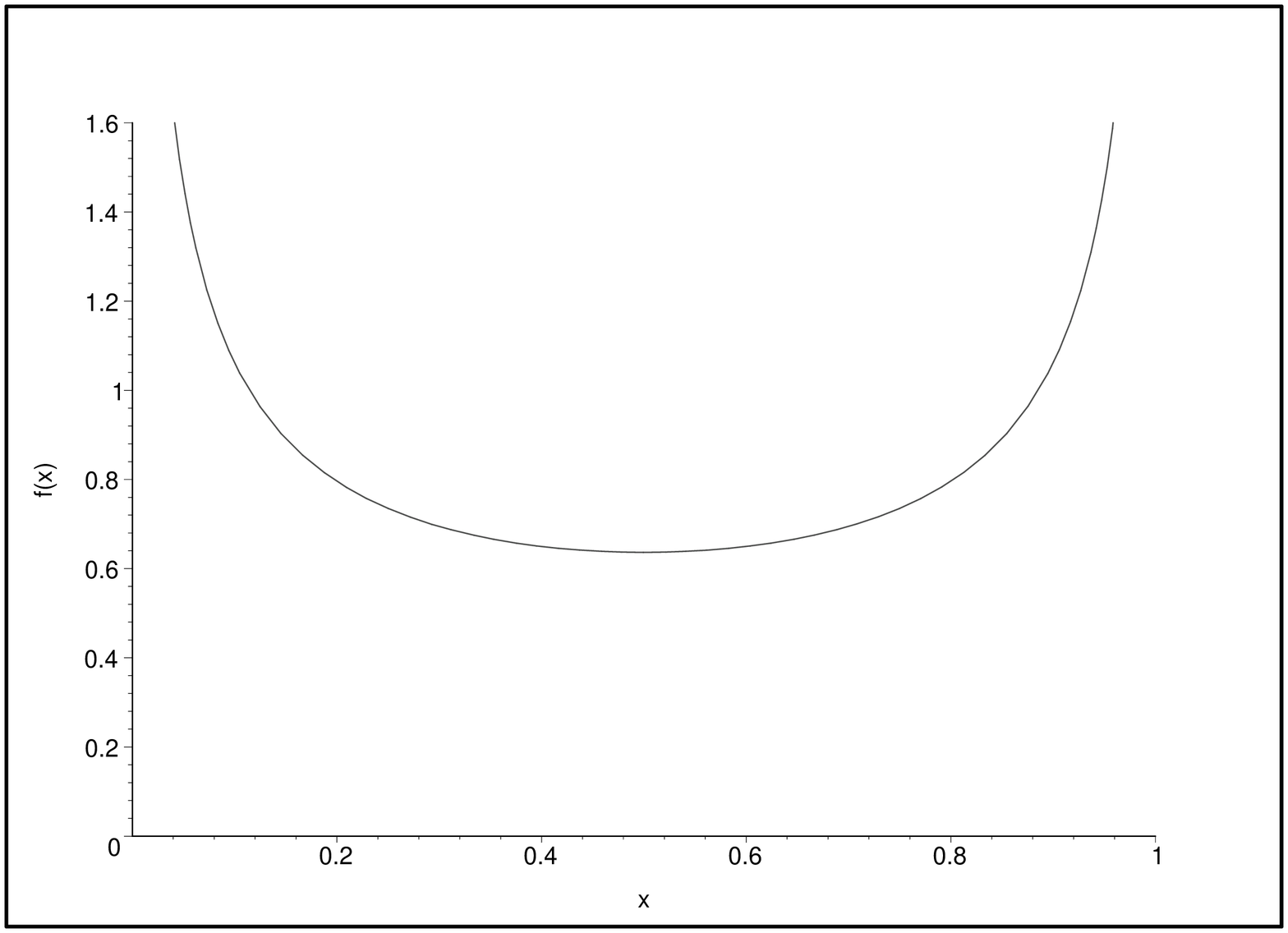, height=180pt , width=140pt}}
\rotatebox{-90}{\epsfig{figure=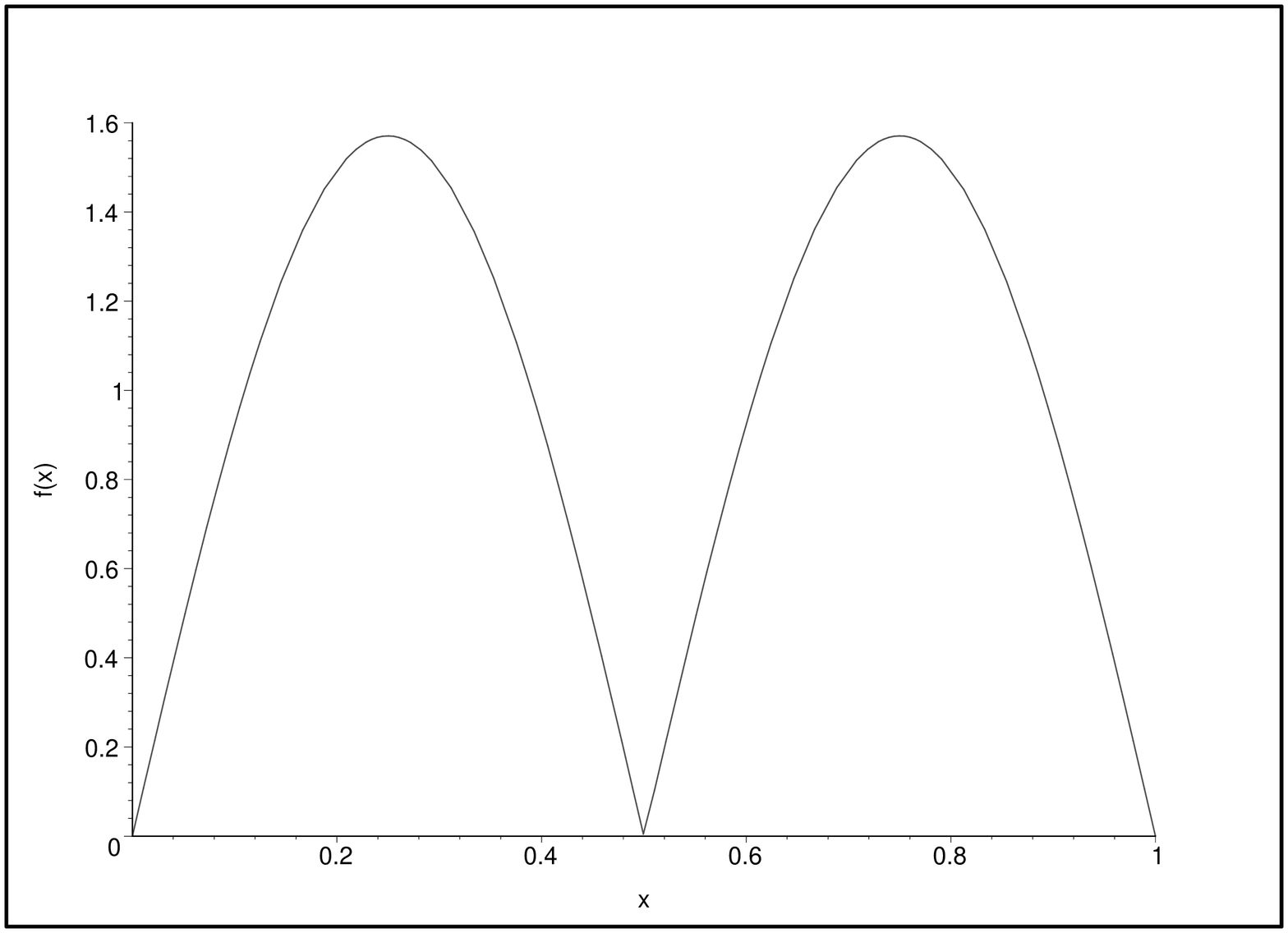, height=180pt , width=140pt}}
\caption{\label{fig:density24-25}
Graph of the density in Examples \ref{ex:arc-sine-density} (left)
and \ref{ex:abs-sine-density} (right).
%Notice that the $y$ axes are differently scaled.
}
\end{figure}

\begin{example}
\label{ex:arc-sine-density}
Consider the distribution with density function
$f(x)=\frac{1}{\pi \sqrt{x\,(1-x)}} \;\I(0<x<1).$
Note that $\Y_2=\{0,1\}$ and $f(x)$ is unbounded at $x \in \{0,1\}$.
See Figure \ref{fig:density24-25} (left) for the plot of $f(x)$.
Instead of $f(x)$, we consider
$g(x)=\frac{\pi\,f(x)}{2\,\arcsin(1-2\delta)}\;\I(\delta<x<1-\delta)$
with cdf $G(x)$.
For $g(x)$, we have $k=\ell=0$ in Theorem \ref{thm:asy-general-Dnm} and then
$\lim_{n \rightarrow \infty}p_n(F)=
\lim_{\delta \rightarrow 0^+}\lim_{n \rightarrow \infty}p_n(G)=1$
using Remark \ref{rem:unbounded}.
The numerically computed value of $p_{1000}(F)$ is $\widehat p_{1000}(F)\approx 1.000$.
$\square$
\end{example}

\begin{remark}
The rate of convergence in  Theorem \ref{thm:kth-order-gen} depends on $f$.
From the proof of Theorem \ref{thm:kth-order-gen},
it follows that for sufficiently large $n$,
$$p_n(F) \approx \frac{f^{(k)}(\y_1^+)\,f^{(\ell)}(\y_2^-)}{\left[f^{(k)}(\y_1^+)+2^{-(k+1)}\,f^{(k)} \left( \frac{\y_1+\y_2)}{2}^+ \right)\right]\,\left[f^{(\ell)}(\y_2^-)+2^{-(\ell+1)}\,f^{(\ell)} \left( \frac{(\y_1+\y_2)}{2}^- \right)\right]} +\frac{c(f)}{n^m},$$
where
$$c(f)=\frac{s_1\,s_3^{\frac{1}{k+1}}\,\G\left( \frac{\ell+2}{\ell+1} \right)+s_2\,s_4^{\frac{1}{\ell+1}}\,\G \left(\frac{k+2}{k+1} \right)}{(k+1)\,(\ell+1)\,s_3^{\frac{k+2}{k+1}}\,s_4^{\frac{\ell+2}{\ell+1}}}$$
with
$\G(x)=\int_{0}^{\infty} \exp(-t)t^{(x-1)}\, dt$ and
\begin{align*}
s_1&=\frac{1}{n^{\frac{k+\ell+1}{\ell+1}}}\,
\frac{(-1)^{\ell+1}}{k!\,(\ell+1)!}\,f^{(k)}(\y_1^+)\,f^{(\ell+1)}(\y_2^-),&
s_3&=\frac{1}{(k+1)!}\left( f^{(k)}(\y_1^+)+2^{-(k+1)}\,f^{(k)}
\left( \frac{(\y_1+\y_2)}{2}^+ \right) \right),\\
s_2&=\frac{1}{n^{\frac{k+\ell+1}{k+1}}}\,\frac{(-1)^\ell}{l!\,(k+1)!}\,
f^{(k+1)}(\y_1^+)\,f^{(\ell)}(\y_2^-),
& s_4&=\frac{(-1)^{\ell+1}}{(\ell+1)!}\left( f^{(\ell)}(\y_2^-)+2^{-(\ell+1)}\,
f^{(\ell)} \left( \frac{(\y_1+\y_2)}{2}^- \right) \right),
\end{align*}
provided the derivatives exist. $\square$
\end{remark}

\begin{figure}[ht]
\centering
\psfrag{density}{}
\psfrag{#}{}
\psfrag{22}{}
\psfrag{23}{}
\psfrag{x}{\scriptsize{$x$}}
\psfrag{f(x)}{\scriptsize{$f(x)$}}
\rotatebox{-90}{\epsfig{figure=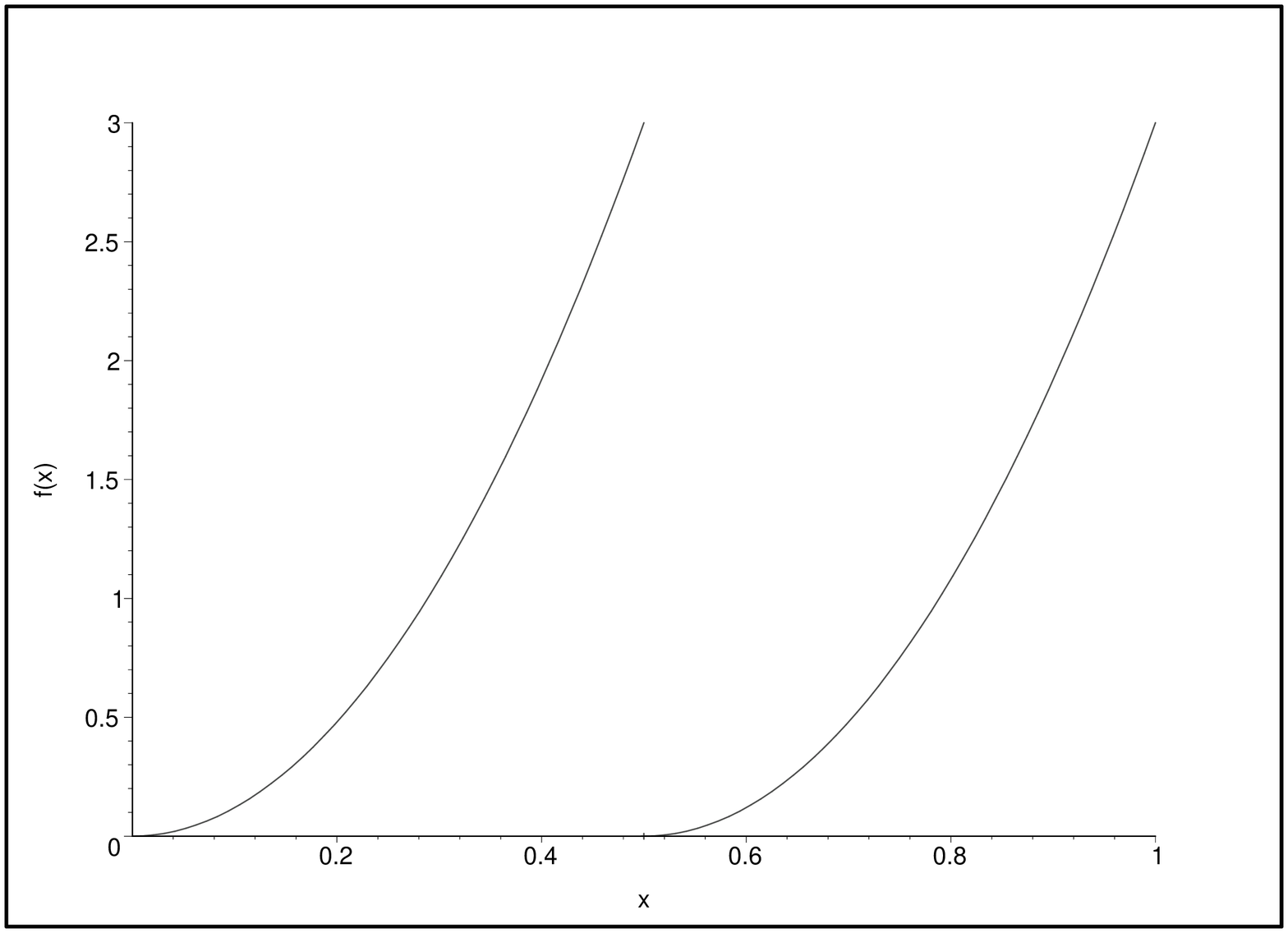, height=180pt , width=140pt}}
\rotatebox{-90}{\epsfig{figure=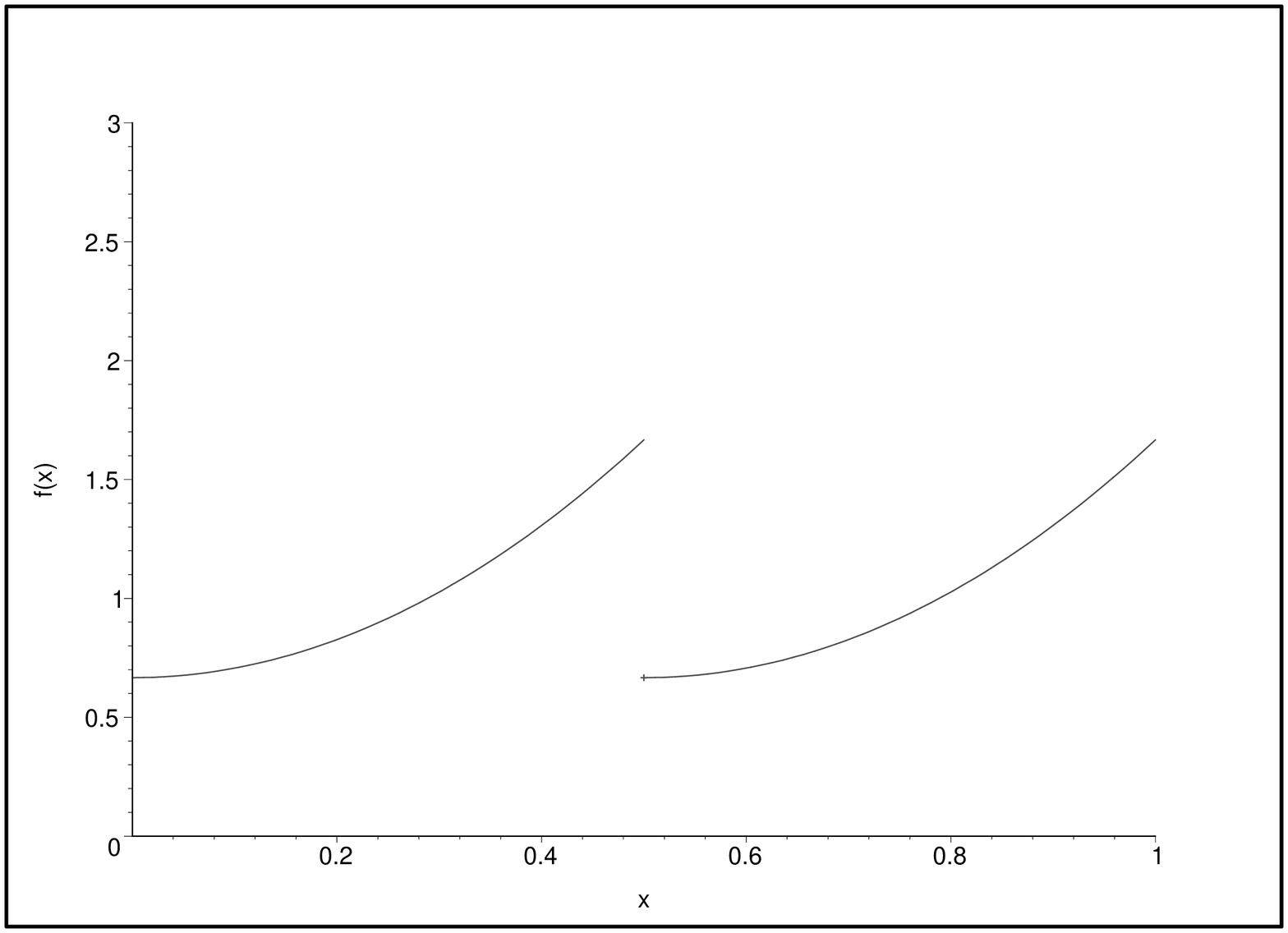, height=180pt , width=140pt}}
\caption{\label{fig:density22-23}
Left plot is for the density in Example \ref{ex:k=q,l=0} with $q=2$
or for the density in Example \ref{ex:PWsmooth3-4} with $\delta=0$.
Right plot is for the density in Example \ref{ex:PWsmooth3-4} with $\delta=2/3$.
}
\end{figure}

\begin{example}
\label{ex:abs-sine-density}
Consider the distribution with absolute sine density
$f(x)=\pi/2 \, | \sin(2\,\pi \,x) | \;\I(0<x<1).$
See Figure \ref{fig:density22-23} (right) for the plot of $f(x)$.
Then $\Y_2=\{0,1\}$ and since
$f(0^+)=f\left( \frac{1}{2}^+ \right)=0$ and $f(1^-)=f\left( \frac{1}{2}^- \right)=0$ and
$f'(0)=f'\left( \frac{1}{2}^+ \right)=\pi^2$ and $f'(1^-)=f'\left( \frac{1}{2}^- \right)=-\pi^2$,
we apply Theorem \ref{thm:kth-order-gen} with $k=\ell=1$.
Then $\lim_{n \rightarrow \infty}p_n(F) =16/25=0.64$.
The numerically computed value (by numerical integration)
of $p_{1000}(F)$ is $\widehat p_{1000}(F) \approx 0.6400$.
$\square$
\end{example}

%\begin{remark}
%Suppose $f^{(k)}(\cdot)$ or $f^{(\ell)}(\cdot)$ is not bounded on $(0,1)$, in particular around $0,1/2,1$, e.g.,
%$\lim_{x \rightarrow 1}f^{(k)}(x)=\infty$.
%Then
%$$ \lim_{n \rightarrow \infty} P\bigl( \g_n(F)=2 \bigr) =\lim_{\delta \rightarrow 0^+}\frac{f^{(k)}(\delta)\,f^{(\ell)}(1-%\delta)}{\bigl(f^{(k)}(\delta)+1/2\,f^{(k)}\bigl( 1/2+\delta \bigr)\bigr)\,\bigl(f^{(\ell)}(1-\delta)+1/2\,f^{(\ell)}\bigl( 1/2-%\delta \bigr)\bigr)}. \;\; \square$$
%\end{remark}

%Notice that
%\begin{itemize}
%\item if $\left( f^{(k)}(0^+),f^{(\ell)}(1^-) \right)=(0,0)$ and
%$f^{(k)} \left( \frac{1}{2}^+ \right)\not=0$ and $f^{(\ell)} \left( \frac{1}{2}^- \right)\not=0$ then $P\bigl(\g_n(F)=2\bigr)%\rightarrow 0$ as $n \rightarrow \infty$, and
%\item  if $f^{(k)}(0^+)\not=0 $ and $f^{(\ell)}(1^-)\not=0$ and $f^{(k)} \left( \frac{1}{2}^+ \right)=f^{(\ell)} \left( \frac{1}%{2}^- \right)=0$ then $P\bigl(\g_n(F)=2\bigr)\rightarrow 1$ as $n \rightarrow \infty$, whose rates are as in the Proof of %Theorem \ref{thm:kth-order}.
%\end{itemize}

The distribution of $\g_n(F)$ depends on the distribution
of $r(X_i)=\min(d(X_i,\y_1),\,d(X_i,\y_2))$.
Based on this, we have the following symmetry result.
\begin{proposition}
\label{prop:symmetric-dist-1}
Let $F_1$ and $F_2$ be two distributions with support $\mS(F_j)\subseteq (\y_1,\y_2)$ for $j=1,2$
such that $F_1(\y_1+x)=1-F_2\,(\y_2-x)$ for all $x \in (0,\y_2-\y_1)$
(hence $f_1(\y_1+x)=f_2\,(\y_2-x)$).
Also, let $\X^j_n$ be a set of iid random variables from $F_j$ for $j=1,2$.
Then the distributions of $\g_n(F_j)$ are identical for $j=1,2$.
\end{proposition}
{\bfseries Proof:}
%Since we have $F_1(\y_1+x)=1-F_2\,(\y_2-x)$, it follows that $f_1(\y_1+x)=f_2\,(\y_2-x)$
%where $f_j$ is the pdf of $F_j$ for $j=1,2$.
By the change of variable $X=\varphi(U)=\y_2-\y_1-U$ for $U \in (0,\y_2-\y_1)$,
we get $F_2\,(\y_1+u)=1-F_1(\y_2-u)$.
Furthermore, $\varphi(u)$ transforms $\G_1(\X^1_n,\NY)$
into $\G_1(\X^2_n,\NY)$ for $\X^2_n$,
so $P(\g_n(F_j)=2)$ are same for both $j=1,2$.
Hence the desired result follows.
$\blacksquare$

Below are asymptotic distributions of $\g_n(F)$ for various families of distributions.
Recall that $p_F=\lim_{n \rightarrow \infty}p_n(F)=\lim_{n \rightarrow \infty}P\bigl(\g_n(F)=2\bigr)$.
The asymptotic distribution of $\g_n(F)$ is $1+\Bernoulli\bigr(p_F\bigl)$.
For the piecewise constant functions in Section \ref{sec:piecewise-constant},
Theorem \ref{thm:kth-order-gen} applies.
See Section 6.1 in \cite{ceyhan:2004c}.

\begin{example}
\label{ex:ax+b}
Consider the distribution $F$ with density $f(\cdot)$ which is of the form
$$f(x)=(a\,x+b)\,\I\bigl( x \in (0,1) \bigr) \text{ with }|a| \le 2,\; b=1-a/2.$$
So $k=\ell=0$ and $f(0^+)=b$, $f(1^-)=a+b$ and $f\left( \frac{1}{2}^+ \right)=f\left( \frac{1}{2}^- \right)=a/2+b$.
Then by Theorem \ref{thm:kth-order-gen}, we have
$$\lim_{n \rightarrow \infty} p_n(F)= \frac{4-a^2}{9-a^2}=:p_F(a).$$
Note that $p_F(a) \in [0,4/9]$ is continuous in $a$ and decreases as $|a|$ increases.
If $a=0$, then $F=\U(0,1)$, and $p_F(a=0)=4/9$.
Moreover, $p_F(a=\pm2)=0$;
that is, for $a=\pm 2$, the asymptotic distribution of $\g_n(F)$ is degenerate.
$\square$
\end{example}

\begin{example}
\label{ex:scaled-normal}
Consider the normal distribution $\N(\mu,\sigma^2)$ restricted to the interval
$(0,1)$ with $\mu \in \R$ and $\sigma>0$.
Then the corresponding density function is given by
$$f(x,\mu,\sigma)=\kappa\left(\frac{1}{\sqrt{2\pi}\,\sigma}\right)
\exp\left(-\frac{(x-\mu)^2}{2\,\sigma^2}\right) \;\I(0<x<1), $$
where $\kappa =\left[\Phi\left(\frac{1-\mu}{\sigma}\right)-
\Phi\left(\frac{-\mu}{\sigma}\right)\right]^{-1}$
with $\Phi(\cdot)$ being the cdf of the standard normal distribution $\N(0,1)$.
Note that $k=\ell=0$, then by Theorem \ref{thm:kth-order-gen}
$$\lim_{n \rightarrow \infty} p_n(F)=\frac{4}
{\left(2+\exp\left(\frac{4\,\mu-1}{8\,\sigma^2}\right)\right)
\left(2+\exp\left(\frac{3-4\,\mu}{8\,\sigma^2}\right)\right)}=:p_F(\mu,\sigma).$$
Observe that $p_F(\mu,\sigma) \in [0,4/9)$ is continuous in $\mu$
and $\sigma$ and increases as $\sigma$ increases for fixed $\mu$.
Furthermore, for fixed $\mu$, $\lim_{\sigma \rightarrow
\infty}p_F(\mu,\sigma)=4/9$ and $\lim_{\sigma \rightarrow
0}p_F(\mu,\sigma)=0$. For fixed $\sigma>0$, $\lim_{\mu \rightarrow
\pm \infty}p_F(\mu,\sigma)=0$, $p_F(\mu,\sigma)$ decreases as
$|\mu-1/2|$ increases, and $p_F(\mu,\sigma)$ is maximized at
$\mu=1/2$. $\square$
\end{example}

\begin{example}
\label{ex:k=q,l=0}
Consider the distribution $F$ with density $f(\cdot)$ which is of the form
$$f(x)=2^q(q+1)\,\left[ x^q\,\I\bigl( 0<x<1/2 \bigr)+(x-1/2)^q\,\I\bigl( 1/2 \le x<1 \bigr) \right]
\text{ with } q \in [0,\infty].$$
%That is, $f_1(x)=2^q(q+1)\,x^q$ and $f_2\,(x)=2^q(q+1)\,(x-1/2)^q$.
See Figure \ref{fig:density22-23} (left) with $q=2$.
Since $f(0^+)=f\left( \frac{1}{2}^+\right)=0$, we can apply
Theorem \ref{thm:kth-order-gen} with $k=q$ and $l=0$.
Then $f^{(q)}(0^+)=(q+1)!\,2^q$, $f(1^-)=(q+1)$,
$f\left( \frac{1}{2}^-\right)=(q+1)$, and
$f ^{(q)}\,\left( \frac{1}{2}^+ \right)=(q+1)!\,2^q$.
By Theorem \ref{thm:kth-order-gen}, we have
$$\lim_{n \rightarrow \infty} p_n(F)= \frac{2^{q+2}}{3\,(1+2^{q+1})}=:p_F(q).$$
Note that $p_F(q) \in [4/9,2/3]$ is a continuous increasing function of $q$.
If $q=0$, then $F=\U(0,1)$.
$\square$
\end{example}

\begin{example}
\label{ex:PWsmooth3-4}
Consider the distribution $F$ with density $f(\cdot)$ which is of the form
$$f(x)=\left( \delta+12\,(1-\delta)\,x^2 \right)\,\I\bigl( 0<x<1/2 \bigr)+
\left( \delta+12\,(1-\delta)\,(x-1/2)^2 \right)\,\I\bigl( 1/2 \le x<1 \bigr)
\text{ with } \delta \in [0,1].$$
%That is,
%$f_1(x)=\bigl( \delta+12\,(1-\delta)\,x^2 \bigr)$ and
%$f_2\,(x)=\bigl( \delta+12\,(1-\delta)\,(x-1/2)^2 \bigr)$.
See Figure \ref{fig:density22-23} with $\delta=0$ (left) and $\delta=2/3$ (right).
Since $f(0^+)=\delta$, $f\,(1^-)=(3-2\,\delta)$,
$f\left( \frac{1}{2}^-\right)=(3-2\,\delta)$, and
$f\left( \frac{1}{2}^+\right)=\delta$, for $\delta \in (0,1]\}$
we have $k=\ell=0$ and so
by Theorem \ref{thm:kth-order-gen}
$$\lim_{n \rightarrow \infty} p_n(F)= 4/9 \text{ for $\delta \in (0,1]$}.$$
Note that if $\delta=1$, then $F=\U(0,1)$.
% but the limiting distribution $\g(F)$ for all $\delta \in [-1,1]\setminus \{0\}$
%is same as that of the $g(\U(0,1))$.
For $\delta= 0$, we can apply Theorem \ref{thm:kth-order-gen} with $k=2$ and $l=0$.
Hence we get $p_F(\delta=0)=16/27$.
Observe that in this example,
$\g_n(F)$ has two distinct non-degenerate distributions at different values of $\delta$.
$\square$
\end{example}

\begin{remark}
\label{rem:same-asydist-uniform}
If, in Theorem \ref{thm:kth-order-gen}, we have $f^{(k)}(0^+)=f^{(k)} \left( \frac{1}{2}^+ \right)$ and $f^{(\ell)}(1^-)=f^{(\ell)} \left( \frac{1}{2}^- \right)$, then
$$\lim_{n \rightarrow \infty}p_n(F)=\left(\frac{1}{1+2^{-(k+1)}}\right)\left(\frac{1}{1+2^{-(\ell+1)}}\right).$$
In particular, if $k=\ell=0$, then
$\lim_{n \rightarrow \infty}p_n(F)=4/9$
(i.e., $\g_n(F)$ and $\g_n\bigl( \U(0,1) \bigr)$ have the same asymptotic distributions).
$\square$
\end{remark}

\begin{example}
\label{ex:beta}
$\text{Beta}(\nu_1,\nu_2)$ with $\nu_1, \nu_2 \ge 1$.  The density function is
$$f(x,{\nu_1},{\nu_2})= {\frac {x^{\nu_1-1}(1-x)^{\nu_2-1}}{\beta(\nu_1,\nu_2) }} \;\I(0<x<1) \text{ where }\beta(\nu_1,\nu_2)=\frac{\Gamma(\nu_1)\,\Gamma(\nu_2)}{ \Gamma(\nu_1 +\nu_2)}.$$
 Then
$\lim_{n \rightarrow \infty}p_n(\text{Beta}(\nu_1,\nu_2))=0$ at rate $O\bigl( n^{-(\nu_1+\nu_2-2)} \bigr)$.
Let $p_n(\nu_1,\nu_2)$ denote the $P\bigl( \g_n(F)=2 \bigr)$ for $F=\text{Beta}(\nu_1,\nu_2)$.
The numerically computed values of $p_n(\nu_1,\nu_2)$ for
$n=1000$ are $\widehat p_{1000}(4,1)=\widehat p_{1000}(1,4) \approx 0.000005$,
$\widehat p_{1000}(4,2)=\widehat p_{1000}(2,4) < 0.00001$ and
$\widehat p_{1000}(2,2) \approx 0.000001$.
$\square$
\end{example}

Here is an example with general support $(\y_1,\y_2)$.
\begin{example}
Consider the distribution $F$ with density $f(\cdot)$ which is of the form
$f(x)=a\,x+b \text{ with } b=\frac{1}{(\y_2-\y_1)}\bigl( 1-a\,(\y_2^2-\y_1^2)/2 \bigr) \text{ and }
|a|\le \frac{2}{(\y_2-\y_1)^2}.$
Using Theorem \ref{thm:kth-order-gen},
we obtain $p_F=\frac{a^2\,(\y_2-\y_1)^4-4}{a^2\,(\y_2-\y_1)^4-9}$.
If $(\y_1,\y_2)=(0,1)$, then $b=1-a/2$ and $p_F(a)=\frac{a^2-4}{a^2-9}$.
In both cases, $p_F(a)$ is maximized for the uniform case; i.e., when $a=0$,
then we have $p_F(a=0)=4/9$.
Furthermore, $\g_n(F)$ is degenerate in the limit when $a=\pm \frac{2}{(\y_2-\y_1)^2}$,
since $p_n(F) \rightarrow 0$ as $n \rightarrow \infty$ at rate $O\bigl( n^{-1} \bigr)$.
%We can also view $a\,x+b$ as a departure from uniformness as
%$\frac{1}{(\y_2-\y_1)} \,\bigl( 2\delta\,x+1-\delta(\y_1+\y_2) \bigr)$ with
%$|\delta| \le 1/(\y_2-\y_1)$.
%Then $p_F=\frac{4\,(\delta^2\,(\y_2-\y_1)^2-1)}{(4\,\delta^2\,(\y_2-\y_1)^2)-9}$
%which equals $4/9$ at $\delta = 0$.
$\square$
\end{example}

For more detail and examples, see Sections 6.4 and 7.1 in \cite{ceyhan:2004c}.

\section{The Distribution of the Domination Number of $\D_{n,m}$-digraphs}
In this section, we attempt the more challenging case of $m>2$.
For $c<d$ in $\R$, define the family of distributions
$$
\mathscr H(\R):=\bigl \{ F_{X,Y}:\;(X_i,Y_i) \sim F_{X,Y} \text{ with support }
\mS(F_{X,Y})=(c,d)^2 \subsetneq \R^2,\;\;X_i \sim F_X \text{ and } Y_i \stackrel{iid}{\sim}F_Y \bigr\}.
$$
We provide the exact distribution of $\g(D_{n,m})$ for
$\mathscr H(\R)$-random digraphs in the following theorem.
Let $[m]:=\bigl\{ 0,1,\ldots,m-1 \bigr\}$ and
$\Theta^S_{a,b}:=\bigl\{ (u_1,\ldots u_b):\;\sum_{i=1}^{b}u_i = a:\; u_i \in S, \;\;\forall i \bigr\}$.
Let $\Y_m=\bigl\{ Y_1,Y_2,\ldots,Y_m \bigl\}$ whose
order statistics are denoted as $Y_{(j)}$ for $j=1,2,\ldots,m$.
Note that the order statistics are distinct a.s.
provided $Y$ has a continuous distribution.
Let $\g(D^j)$ be the domination number of the digraph induced by $\X^j$ and $\Y^j$
(see Section \ref{sec:domination-number-Dnm}).
Given $Y_{(j)}=\y_{(j)}$ for $j=1,\ldots,m$,
let $F_j$ be the (conditional) marginal distribution of $X$ restricted to
$\mI_j=\left(\y_{(j-1)},\y_{(j)}\right)$ for $j=1,\ldots,(m+1)$,
$\vec{n}$ be the vector of numbers of $\X$ points $n_j$
falling into intervals $\mI_j$.
%To make the dependence on $F_j$ explicit we denote $P(\g(D^j)=2)$ as
%$p_{n_j}(F_j)$ and $\lim_{n_j \rightarrow \infty}p_{n_j}(F_j)=:p_{F_j}$.
Let $f_{\vec{Y}}(\vec{\y})$ be the joint distribution of the order statistics of $\Y_m$,
i.e., $f_{\vec{Y}}(\vec{\y})=\frac{1}{m!}\prod_{j=1}^m f(\y_j)\,\I(c<\y_1<\ldots<\y_m<d)$,
and $f_{j,k}(\y_j,\y_k)$ be the joint distribution of $Y_{(j)},Y_{(k)}$.
Then we have the following theorem which is a generalization of the main result of \cite{priebe:2001}.

\begin{theorem}
\label{thm:general-Dnm}
Let $D$ be an $\mathscr H(\R)$-random $\D_{n,m}$-digraph.
Then the probability mass function of the domination number of D is given by
$$P(\g(D_{n,m})=k)=\int_{\mathscr S} \sum_{\vec{n} \in \Theta^{[n+1]}_{n,(m+1)}}
\sum_{\vec{k}\in \Theta^{[3]}_{k,(m+1)}} P(\vec{N}=\vec{n})\,\zeta(k_1,n_1)\,\zeta(k_{m+1},\,n_{m+1})
\prod_{j=2}^{m}\eta(k_j,n_j)f_{\vec{Y}}(\vec{\y})\,d\y_1 \ldots d\y_m$$
where
$P(\vec{N}=\vec{n})$ is the joint probability of $n_j$ points
falling into intervals $\mI_j$ for $j=1,2,\ldots,(m+1)$,
$k_j \in \{0,1,2\}$, and
\begin{align*}
\zeta(k_j,n_j)&=\max\bigl( \I(n_j=k_j=0),\I(n_j \ge k_j=1) \bigr) \text{ for } j=1,(m+1),
\text{ and }\\
\eta(k_j,n_j)&=\max \bigl( \I(n_j=k_j=0),\I(n_j \ge k_j \ge 1) \bigr)\cdot
p_{n_j}(F_j)^{\I(k_j=2)}\,\bigl( 1-p_{n_j}(F_j) \bigr)^{\I(k_j=1)}\\
& \text{ for $j=2,\ldots,m,$ and the region of integration is given by}\\
\mathscr S:=\bigl\{&(\y_1,\y_2,\ldots,\y_m)\in (c,d)^2:\,c<\y_1<\y_2<\ldots<\y_m<d \bigr\}.
\end{align*}
\end{theorem}
{\bfseries Proof:}
For $\g(D_{n,m})=\sum_{j=1}^{m+1}\,\g(D^j)=k$,
we must have $\g(D^j)=k_j$ for $j=1,\ldots,(m+1)$ so that
$\sum_{j=1}^{m+1}\,k_j=k$ and $\sum_{j=1}^{m+1}n_j=n$.
By definition, $\Theta^{[n+1]}_{n,(m+1)}$ is the collection
of such $\vec{n}$ and since $k_j \in \{0,1,2\}$
for all $j=1,\ldots,(m+1)$, $\Theta^{[3]}_{k,(m+1)}$ is the
collection of such $\vec{k}$.
We treat the end intervals, $\mI_1$ and $\mI_{m+1}$, separately.
The indicator functions in the statement of the theorem guarantees
that the pairs $n_j,k_j$ are compatible for $j\in \{1,(m+1)\}$;
that is, incompatible pairs such as $(n_j=0,k_j>0)$ are eliminated.
The $\zeta$ terms equal unity if $(n_j,k_j)$ are compatible.
Therefore we have
\begin{multline*}
P(\g(D_{n,m})=k)=\int_{\mathscr S} \sum_{\vec{n} \in \Theta^{[n+1]}_{n,(m+1)}}
\sum_{\vec{k}\in \Theta^{[3]}_{k,(m+1)}} P(\vec{N}=\vec{n})\,
\prod_{j=1}^{m+1}\eta(k_j,n_j)f_{\vec{Y}}(\vec{\y})\,d\y_1 \ldots d\y_m\\
= \int_{\mathscr S} \sum_{\vec{n} \in \Theta^{[n+1]}_{n,(m+1)}}
\sum_{\vec{k}\in \Theta^{[3]}_{k,(m+1)}} P(\vec{N}=\vec{n})\,\prod_{j\in\{1,(m+1)\}}\eta(k_j,n_j)
\prod_{j=2}^{m}\eta(k_j,n_j)f_{\vec{Y}}(\vec{\y})\,d\y_1 \ldots d\y_m \\
=\int_{\mathscr S} \sum_{\vec{n} \in \Theta^{[n+1]}_{n,(m+1)}}
\sum_{\vec{k}\in \Theta^{[3]}_{k,(m+1)}} P(\vec{N}=\vec{n})\,\zeta(k_1,n_1)\,
\zeta(k_{m+1},\,n_{m+1})
\prod_{j=2}^{m}\eta(k_j,n_j)f_{\vec{Y}}(\vec{\y})\,d\y_1 \ldots d\y_m
\end{multline*}
where we have used the conditional pairwise independence of $\g(D^j)$.
The $\eta$ terms are based on the compatibility of pairs $(n_j,k_j)$ for
$j=1,\ldots,(m+1)$ and $p_{n_j}(F_j)=P(\g(D^j)=2)$.
$\blacksquare$

For $n,m < \infty$, the expected value of domination number is
\begin{equation}
\label{eqn:expected-gamma-Dnm}
\E[\g(D_{n,m})]=P\left(X_{(1)}<Y_{(1)}\right)+P\left(X_{(n)} > Y_{(m)}\right)+
\sum_{j=2}^m\sum_{k=1}^n\,P(N_j=k)\,\E[\g(D^j)]
\end{equation}
where
\begin{multline*}
P(N_j=k)=\\
\int_c^d\int_{\y_{(j-1)}}^d
f_{j-1,j}\left(\y_{(j-1)},\y_{(j)}\right) \Bigl[F_X\left(\y_{(j)}\right)-
F_X\left(\y_{(j-1)}\right)\Bigr]^k\Bigl[1-\left(F_X\left(\y_{(j)}\right)-
F_X\left(\y_{(j-1)}\right)\right)\Bigr]^{n-k}\,d\y_{(j)}d\y_{(j-1)}
\end{multline*}
and $\E[\g(D^j)]=1+p_{k}(F_j)$.

\begin{corollary}
\label{cor:Egnn goes infty}
For $F_{X,Y} \in \mathscr H(\R)$ with support $\mS(F_X) \cap \mS(F_Y)$
of positive measure, $\lim_{n \rightarrow \infty}\E[\g(D_{n,n})] = \infty$.
\end{corollary}
{\bfseries Proof:}
Consider the intersection of the supports $\mS(F_X) \cap \mS(F_Y)$
that has positive (Lebesgue) measure.
For $\mS(Y) \setminus \mS(X)$; i.e., in the intervals $\mI_j$ falling outside
the intersection $\mS(F_X) \cap \mS(F_Y)$, the domination number of the component
$D^j$ is $\g(D^j)=0$  w.p. 1 but inside the intersection,
$\g(D^j)>0$ w.p. 1 for infinitely many $j$.
That is,
\begin{eqnarray*}
\E[\g(D_{n,n})]&=&P\left(X_{(1)}<Y_{(1)}\right)+P(X_{(n)}>Y_{(n)})+
\sum_{j=2}^n\sum_{k=1}^n\,P(N_j=k)\,\E[\g_{N_j}(F_j)]\\
&>&\sum_{j=2}^n\sum_{k=1}^n\,P(N_j=k)\,\E[\g_{N_j}(F_j)]=
\sum_{j=2}^n\sum_{k=1}^n\,P(N_j=k)\,(1+p_{N_j}(F_j))\\
&>&\sum_{j=2}^n\sum_{k=1}^n\,P(N_j=k) > \sum_{j=2}^n\,P(N_j \ge 1)\\
&\approx& n \text{~~~ (for sufficiently large $n$)}
\end{eqnarray*}
where $\E[\g_{N_j}(F_j)]=(1+p_{N_j}(F_j))$ follows from the fact that
$\g_{N_j}(F_j)\sim 1+\Bernoulli(p_{N_j}(F_j))$ from Theorem \ref{thm:gamma 1 or 2}.
Furthermore, $P(N_j \ge 1) \approx 1$ for sufficiently large $n$.
Then the desired result follows.
$\blacksquare$

\begin{theorem}
\label{thm:asy-general-Dnm}
Let $D_{n,m}$ be an $\mathscr H(\R)$-random $\D_{n,m}$-digraph.
Then (i) for fixed $n<\infty$, $\lim_{m \rightarrow \infty}\g(D_{n,m})=n$ a.s.
(ii) for fixed $m<\infty$,
$\lim_{n \rightarrow \infty}\g(D_{n,m})\stackrel{d}{=}m+1+\sum_{j=1}^m B_j$,
where $B_j\sim \Bernoulli(p_{F_j})$ where $\stackrel{d}{=}$ stands for equality in distribution.
\end{theorem}
{\bfseries Proof:} Part (i) is trivial.
As for part (ii), first note that $N_j \rightarrow \infty$ as $n \rightarrow \infty$ for all $j$ a.s.,
hence $\lim_{n\rightarrow \infty} \g(D^1)=\lim_{n\rightarrow \infty} \g(D^{m+1})=1$ a.s.
and $\lim_{n\rightarrow \infty} \g(D^j)=1+\Bernoulli(p_{F_j})$ a.s. for $j=2,\ldots,m$ where
$$p_{F_j}=\int_c^d\int_{\y_{(j-1)}}^d H^*\left(\y_{(j-1)},\y_{(j)}\right)\,f_{j-1,j}\left(\y_{(j-1)},\y_{(j)}\right)\,d\y_{(j)}d\y_{(j-1)}$$
with $H^*\left(\y_{(j-1)},\y_{(j)}\right)=\lim_{n_j \rightarrow \infty}(p_{n_j}(F_j))$
which is given in Theorem \ref{thm:general-Dnm} for $F_j$
with density $f_j$ whose support is $\left( \y_{(j-1)},\y_{(j)}\right)$.
Then the desired result follows.
$\blacksquare$

So far, $\Y_m$ is assumed to be a random sample,
so $P(\g(D_{n,m})=k)$ includes the integration with respect to $f_{\vec{Y}}(\vec{\y})$
which can be lifted by conditioning.
%Next we do the calculations given $\Y_m$.
%\begin{corollary}
Conditional on $\Y_m=\left \{\y_{(1)},\ldots,\y_{(m)} \right \}$,
by Theorem \ref{thm:general-Dnm} we have
$$P(\g(D_{n,m})=k)= \sum_{\vec{n} \in \Theta^{[n+1]}_{n,(m+1)}}
\sum_{\vec{k}\in \Theta^{[3]}_{k,(m+1)}} P(\vec{N}=\vec{n})\,\zeta(k_1,n_1)\,\zeta(k_{m+1},\,n_{m+1})
\prod_{j=2}^{m}\eta(k_j,n_j),$$
where $\zeta(k_j,n_j)$ and $\eta(k_j,n_j)$ are  as in Theorem \ref{thm:general-Dnm};
%\end{corollary}
%{\bfseries Proof:} Similar to the Proof of Theorem \ref{thm:general-Dnm}.
%$\blacksquare$
%Conditional on $\Y_m=\left \{\y_{(1)},\ldots,\y_{(m)} \right \}$,
and the expected domination number $\E[\g(D_{n,m})]$ is as in
Equation \eqref{eqn:expected-gamma-Dnm} with
$P(N_j=k)=\left[F_X\left(\y_{(j)}\right)-F_X\left(\y_{(j-1)}\right)\right]^k
\left[1-\left(F_X\left(\y_{(j)}\right)-F_X\left(\y_{(j-1)}\right)\right)\right]^{n-k}$;
and
%\begin{theorem}
%Let $D_{n,m}$ be a $\mathscr H(\R)$-random $\D_{n,m}$-digraph.
%Then conditional on $\Y_m=\left \{\y_{(1)},\ldots,\y_{(m)} \right \}$,
$\lim_{n \rightarrow \infty}\g(D_{n,m})\stackrel{d}{=}m+1+\sum_{j=1}^m B_j$,
where $B_j\sim \Bernoulli(p_{F_j})$ with $p_{F_j}:=\lim_{n_j \rightarrow \infty}p_{n_j}(F_j)$.
%\end{theorem}
%{\bfseries Proof:} Similar to part (ii) of Proof of Theorem \ref{thm:asy-general-Dnm},
%with $p_{F_j}=\lim_{n_j \rightarrow \infty}(p_{n_j}(F_j))$ which is given in
%Theorem \ref{thm:kth-order-gen} for $F_j$ with density $f_j$
%whose support is $\left(\y_{(j-1)},\y_{(j)}\right)$.
%$\blacksquare$

Let $F_X$ be a distribution with support $\mS(F_X)\subseteq(0,1)$ and density $f_X(x)$.
Conditional on $\Y_m=\left \{\y_{(1)},\ldots,\y_{(m)} \right \}$,
let $F_j$ be the distribution with density
$f_j(x)=\frac{1}{\left(\y_{(j)}-\y_{(j-1)}\right)}
f_X\left( \frac{x-\y_{(j-1)}}{\y_{(j)}-\y_{(j-1)}} \right)$ for $j=2,\ldots,m$,
and $\mS(F_j(x))$ is non-empty for $j\in\{1,(m+1)\}$.
By this construction, the independence of the distribution of
$\g_n(F_j)$ from $\mI_j$ obtains; i.e.,
$\g_n(F_j)\stackrel{d}{=}\g_n(F_X)$ for all $j \in \{1,\ldots,(m+1)\}$.
Now consider the family $\mathscr H_\U(\R)$ defined as
$$
\mathscr H_\U(\R):=\bigl\{ F_{X,Y}:(X_i,Y_i) \sim F_{X,Y},\;Y_j \stackrel{iid}{\sim}\U(c,d)
\text{ for }(c,d) \subsetneq \R, \text{ and } X_i|\Y_m \stackrel{iid}{\sim} F_j \bigr\}.
$$
Clearly $\mathscr H_\U(\R) \subsetneq \mathscr H(\R)$.

\begin{corollary}
\label{cor:XF-YU-Dnm}
Suppose $F_{X,Y} \in \mathscr H_\U(\R)$.  Then
$$P(\g(D_{n,m})=k)= \sum_{\vec{n} \in \Theta^{[n+1]}_{n,(m+1)}}
 \sum_{\vec{k}\in \Theta^{[3]}_{k,(m+1)}} P(\vec{N}=\vec{n})\,\zeta(k_1,n_1)\,\zeta(k_{m+1},\,n_{m+1})
\prod_{j=2}^{m}\eta(k_j,n_j)$$
where $\zeta(k_j,n_j)$ and $\eta(k_j,n_j)$ are as in Theorem \ref{thm:general-Dnm}.
\end{corollary}
Note that if in addition, $P_{F_j}(X \in \mI_j)=P_\U(X \in \mI_j)$ for all $j$,
then $P(\vec{N}=\vec{n})={n+m \choose n}^{-1}$,
since each $\vec{n} \in \Theta^{[n+1]}_{n,(m+1)}$ occurs with probability ${n+m \choose n}^{-1}$.
Moreover, $F=\U(c,d)$ is a special case of Corollary \ref{cor:XF-YU-Dnm}.
For $n,m <\infty$, we have the explicit form of $p_{n_j}(F_j)$ for $F_j$
with piecewise constant density $f_j$.

Here are some examples which are generalized from piecewise-constant densities
so that now the distribution of $\g(D^j)$ is independent from the support $(\y_{(j-1)},\y_{(j)})$.
Hence Corollary \ref{cor:XF-YU-Dnm} applies to these examples.
\begin{example}
\begin{itemize}
\item[]
Let $u_j:=\frac{\bigl(\y_{(j-1)}+\y_{(j)}\bigr)}{2}$ and $v_j:=\y_{(j)}-\y_{(j-1)}$.
\item
If $f(\cdot)$ is of the form
$$f(x)=\frac{1}{(1-2\,\delta)\,v_j}\,
\I \left( x \in \left( \y_{(j-1)}+\delta\,v_j, \y_{(j)}-\delta\,v_j\right) \right)
\text{ with }\delta \in [0,1/3]$$ then $p_n(F)$ is as in Equation \eqref{eqn:pwc1}.

\item
If $f(\cdot)$ is of the form
$$
f(x)=\frac{1}{(1-2\,\delta)\,v_j}\,\I \left( x \in \left( \y_{(j-1)},u_j-\delta\,v_j \right)
\cup \left[ u_j+\delta\,v_j,\y_{(j)}\right) \right) \text{ with }\delta \in [0,1/3],
$$
then $p_n(F)$ is as in Equation \eqref{eqn:pwc2}.

\item
If $f(\cdot)$ is of the form
$$f(x)=\frac{(1+\delta)}{\,v_j}\,
\I\left(x \in \left(\y_{(j-1)},u_j\right)\right)+
\frac{(1-\delta)}{\,v_j}\, \I\left(x \in
\left[u_j,\y_{(j)}\right) \right),$$ then $p_n(F)$ is as in Equation \eqref{eqn:pwc3}.

\item
If $f(\cdot)$ is of the form
$$
f(x)=f_1(x)\,\I\left( x \in \left(\y_{(j-1)},t_j\right)\right)+
f_2\,(x)\,\I\left(x \in \left[t_j, w_j\right)\right)+
f_3\,(x)\,\I\left( x \in \left[w_j,\y_{(j)} \right)\right)
$$
where $t_j=\frac{\y_{(j)}+3\,\y_{(j-1)}}{4}$,
$w_j=\frac{3\y_{(j)}+\y_{(j-1)}}{4}$,
$f_1(x)=\frac{(1+\delta)}{\,v_j}$,
$f_2\,(x)=\frac{(1-\delta)}{\,v_j}$ and
$f_3\,(x)=\frac{(1+\delta)}{\,v_j}$, then $p_n(F)$ is as in Example \ref{ex:PW4}. $\square$
\end{itemize}
\end{example}

\begin{theorem}
Let $D$ be an $\mathscr H_\U(\R)$-random $\D_{n,m}$-digraph with the additional
assumption that $P_{F_j}(X \in \mI_j)=P_\U(X \in \mI_j)$ for all $j$.
Then
$$\E[\g(D_{n,m})]=\frac{2\,n}{n+m}+\frac{n!\,m\,(m-1)}{(n+m)!}\,
\sum_{i=1}^{n}\frac{(n+m-i-1)!}{(n-i)!}\,(1+p_{i}(F))$$
where $p_{i}(F)=P(\g(D_{i,2})=2)$.
\end{theorem}
{\bfseries Proof:} Similar to the Proof of Theorem 4 in \cite{priebe:2001}.
$\blacksquare$

Furthermore, from Corollary \ref{cor:Egnn goes infty},
we have $\E[\g(D_{n,n})]\rightarrow \infty$ as $n\rightarrow \infty$.

\begin{theorem}
 Let $D_{n,m}$ be an $\mathscr H_\U(\R)$-random $\mathscr D_{n,m}$-digraph.
Then (i) for fixed $n<\infty$, $\lim_{m \rightarrow \infty}\g(D_{n,m})=n$ a.s.
(ii) for fixed $m<\infty$, $\lim_{n \rightarrow \infty}\g(D_{n,m})\stackrel{d}{=}m+1+B$,
where $B\sim Binomial(m-1,p_F)$ where $p_F=\lim_{n\rightarrow \infty}P(\g(D_{n,2})=2)$.
\end{theorem}
{\bfseries Proof:} Similar to the Proof of Theorem 5 in \cite{priebe:2001}.
$\blacksquare$

\begin{remark}
\label{rem:ext-multi-dim}
\textbf{Extension to Multi-dimensional Case:}
The existence of ordering of points in $\R$ is crucial in our calculations.
The order statistics of $\Y_m$ partition the support $(c,d)$ into disjoint intervals
a.s. which can also be viewed as the Delaunay tessellation of $\R$ based on $\Y_m$.
This nice structure in $\R$ avails a minimum dominating set and hence the domination
number, both in polynomial time.
Furthermore, the $\G_1$-region is readily available by the order statistics of $\X_n$;
also the components of the digraph restricted to intervals $\mI_j$ (see Section \ref{sec:domination-number-Dnm})
are not connected to each other, since $\NY(x)\cap\NY(y)=\emptyset$
for $x,\,y$ in distinct intervals.
The straightforward extension to multiple dimensions (i.e., $\R^d$ with $d>1$)
does not have a nice ordering structure;
and $\Y_m$ does not readily partition the support,
but we can use the Delaunay tessellation based on $\Y_m$.
Furthermore, in multiple dimensions finding a minimum dominating set is an NP-hard problem;
and $\G_1$-regions are not readily available (in fact for $n_j>3$,
complexity of finding the $\G_1$-regions is an open problem).
In addition, in multiple dimensions the components of the digraph restricted to Delaunay cells
are not necessarily disconnected from each other,
since $\NY(x)\cap\NY(y)\not=\emptyset$
might hold for $x,\,y$ in distinct Delaunay cells.
These have motivated us to generalize the proximity map $\NY$
in order to avoid the difficulties above.
See \cite{ceyhan:2003a, ceyhan:2005e},
where two new families of proximity maps are introduced,
and the generalization of CCCD are called proximity catch digraphs.
The distribution of the domination number of these proximity maps
is still a topic of ongoing research. $\square$
\end{remark}

\section{Discussion}
This article generalizes the main result of \cite{priebe:2001} in several directions.
\cite{priebe:2001} provided the exact (finite sample) distribution of the
class cover catch digraphs (CCCDs) based on $\X_n$ and $\Y_m$
both of which were sets of iid random variables from
a uniform distribution on $(c,d) \subset \R$ with $-\infty<c<d<\infty$
and the proximity map $\NY(x):=B(x,r(x))$ where $r(x):=\min_{\y \in \Y_m}d(x,\y)$.
First, given $\Y_2=\{\y_1,\y_2\} \subset \R$,
we lift the uniformity assumption of $\X_n$ by assuming it
to be from a non-uniform distribution $F$ with
support $\mS(F) \subseteq (\y_1,\y_2)$.
The exact distribution of the domination number of the associated CCCD,
$\g_n(F)$, is calculated for $F$ that has
piecewise constant density $f$ on $(\y_1,\y_2)$.
For more general $F$, the exact distribution is not analytically available
in simple closed form,
so we compute it by numerical integration.
However, the asymptotic distribution of $\g_n(F)$ is tractable, which is the
one of the main results of this article.
Unfortunately, the distribution of $\g_n(F)$ depends on $\Y_2$,
hence the distribution of the domination number of a CCCD, $\g(D_{n,m})$,
for $\X_n$ and $\Y_m$ with $m>2$,  for general $F$
includes integration with respect to order statistics of $\Y_m$.
We provide the conditions that make $\g(D_{n,m})$ independent of $\Y_m$.
As another generalization direction,
we also devise proximity maps depending on $F$ that will yield the
distribution identical to that of $\g_n(\U(\y_1,\y_2))$.
Our set-up is more general than the one given in \cite{priebe:2001}.
The definition of the proximity map is generalized to any probability space
and is only assumed to have a regional relationship to determine the inclusion
of a point in the proximity region.
%Priebe and several co-authors applied CCCDs in pattern classification and obtained
%relatively good performance (\cite{devinney:2002a}, \cite{marchette:2003},
%\cite{priebe:2003a}, and \cite{priebe:2003a}).

The exact (finite sample) distribution of
$\g_n(F)$ characterizes $F$ up to a special type of symmetry
(see Proposition \ref{prop:symmetric-dist-1}).
Furthermore, this article will form the foundation of the generalizations and calculations
for uniform and non-uniform cases in multiple dimensions.
As in \cite{ceyhan:2005e},
we can use the domination number in testing
one-dimensional spatial point patterns
and our results will help make the power comparisons
possible for large families of distributions.

\section*{Acknowledgments}
I would like to thank the anonymous referees, whose constructive
comments and suggestions greatly improved the presentation and flow
of this article.

%\bibliography{References}
\bibliographystyle{apalike}

\end{document}